\documentclass[preprint,12pt]{elsarticle}

\usepackage{graphicx}
\usepackage{amssymb}
\usepackage{amsmath}
\usepackage{bm}
\usepackage[colorlinks,bookmarksopen,bookmarksnumbered,citecolor=red,urlcolor=red]{hyperref}

\usepackage[utf8]{inputenc}

\newtheorem{Remarque}{Remark}
\newtheorem{Def}{Definition}
\newtheorem{prop}{Proposition}
\newtheorem{lem}{Lemma}

\journal{EJFM/B}
\date{November 7, 2020}
\begin{document}

\begin{frontmatter}


\title{On the sound speed in two-fluid mixtures and the implications for CFD model validation}


\author[label1]{Saad Benjelloun}
\author[label1,label2]{Jean-Michel Ghidaglia}
\address[label1]{Mohammed VI Polytechnic University (UM6P), Morocco, Modeling Simulation \& Data Analytics}
\address[label2]{Université Paris-Saclay, ENS Paris-Saclay, CNRS, Centre Borelli, F-91190, Gif-sur-Yvette, France}

\begin{abstract}
Study of the propagation of sound in a single non-ideal fluid originates with Stokes in 1845 and Kirchhoff in 1868. The situation is much more complex in the case of two-fluid flow, both from the physical point of view, as the configuration of the flow matters greatly, and from the analytical point of view. The principle two-fluid models currently in use for CFD are the focus of this article. It is shown that analytical expressions for the speed of sound depend heavily on the chosen model. These sound speed expressions are compared with experimental values. The consequences for CFD models are discussed in the final section of this paper. It is found that numerical models with inaccurate wave speeds lead to incorrect numerical solutions, despite the accuracy of the numerical scheme.
\end{abstract}


\end{frontmatter} 


\section{Introduction}\label{S:1:intro}


Since the Seventies, Computational Fluid Dynamics (CFD) for single-fluid flow has made tremendous progress and today, at least for the engineering community, a wealth of codes are available both in open source and commercial formats, with the possibility of use in industrial applications. This success has resulted from both theoretical advances (physical modeling, numerical methods and algorithms, software engineering, ...) and exponential growth of computational resources.\\

\noindent The situation for multi-fluid flows is more contrasted because physical modeling of multi-fluid and multi-phase flows remains an issue. Derivation of a physical model for multi-fluid flow that explicitly takes account of the interfaces between the fluids is a standard question. However, when these interfaces become complex, it is not possible to derive effective computational models 
and the need for computational resources exceeds what is available by many orders of magnitude.\\

\noindent State of the art multi-fluid flow computation is reminiscent of the situation for the computation of single fluid turbulent flows. Whilst all the information is present in the Navier-Stokes Equations, the number of degrees of freedom necessary to capture the multi-fluid flow at large Reynolds numbers with this physical model is beyond present computational capacities: both memory and CPU time are strong bottlenecks. The solution to this problem has come from turbulence modeling which provides statistical models {\it i.e.} averaged models like {\it e.g.} ensemble-averaged models. The reader is referred to the text book from S.B. Pope \cite{Pope} for a comprehensive account of this matter.\\
\subsection{Averaged models for two-fluid flows}
For multi-fluid flows, a parallel methodology to that of turbulent flows has led to averaged models, see {\it e.g.} the reference book 
 M. Ishii and T. Hibiki \cite{Ishii}. Since multi-fluid flow involves much more complex phenomena than single fluid flow, there are a large number of physical model types and for each single model a large number of modeling terms must be investigated in order to derive a closed model suitable for discretization and simulation. This is inherent to the so-called {\it regimes} in multi-fluid flow. Regimes are spatial configurations of the interfaces between the fluids. For instance for two non-miscible fluids, a liquid and a gas, we can have stratified flows (the two-fluids are separated by a gradual free surface), dispersed flows (either bubbles or droplets), slug flows (large patch of one of the two-fluids), churn flows (transient and very unstructured flows), ... Moreover, in practical situations, the flows are transient and the spatial configurations evolve with time and space.\\

\noindent The story of multi-phase CFD is interesting by itself. In the sixties, a few groups were actively working on building codes in this area, but the task was made challenging by the lack of theoretical work on the subject. The reader is referred to the personal memoir by R. Lyczkowski \cite{Lyczkowski} for an account of this fascinating story.\\

\noindent  Today, there is unfortunately no ``all purpose" model for the simulation of two-fluid flow. However, a large variety of physical models exist\footnote{Some of the main classes of models will be listed in this article, see Section \ref{sec_2-fluid_mod}.}, where each model has a dedicated range of applications. There are at least two issues that must be investigated prior to building a code from a physical model. The first pertains to the physical relevance of the model and the second issue deals with its well-posedness, {\it i.e.} uniqueness of the time evolution problem and stability with respect to the initial and boundary data. For the latter, the reader is referred to the review paper of H. Stewart and B. Wendroff \cite{Stewart}, where, in particular, hyperbolicity of some models for two-fluid flow is studied. The main purpose of this paper is related to the physical relevance of the model and is focused on the following question: {\it does the considered model produce the correct value for the Speed of Sound (SoS)?}\\

\noindent Actually, this question is of great importance for at least two reasons. Physical relevance is a paramount issue when dealing with CFD but also, since most modern numerical methods use so-called characteristic waves ({\it e.g.} exact or approximate Riemann Solvers), it is critical that the sound waves propagate at a physical speed. Propagation of sound waves in multi-phase mixtures is a mature subject from the experimental point of view, see {\it e.g.} Leroy \cite{Leroy} and references within. Hence, this work is concerned with investigating the main physical Eulerian models for two-fluid flows found in the literature and studying analytically the predicted speed of sound values.

\subsection{Content of this work} After introducing the notation used in this paper (Section \ref{sec_notations}), the usual definition for the speed of sound in a fluid is given in Section \ref{sec_sos}. The case of the classical Navier-Stokes model is used to introduce a typical dispersion relation and values for the sound speed and attenuation are obtained. 
Next, in Section \ref{sec142}, a brief literature review on analytical and experimental results for the speed of sound in two-fluid mixtures is presented. Section \ref{sec_2-fluid_mod} lists the main Eulerian models used in the numerical modeling of two-fluid flows. 
The convective part of these models, {\it i.e.} without diffusion, surface tension or change of phase effects, is first dealt with. Mainly, these models can be divided into single-velocity models and two-velocities models. In Section \ref{Sec3} the value of the speed of sound for each of the models from the previous Section is found. In Section \ref{Sec4} the effect of molecular diffusion (viscosity of the fluids) is studied by introducing the dispersion relation for these models. 
In Section \ref{Sec5} the consequences of our results for CFD models are discussed. Finally, after proposing some directions for future work (Section \ref{Sec6}), technical results are collected in the Appendices at the end.

\subsection{Notation}\label{sec_notations}
The notation presented in the following $3$ tables will be used throughout this paper. In Table \ref{tab4}, the notation used for a single fluid model is given.

\begin{table}[ht!]
    \centering
    \begin{tabular}{|c|l||c|l|}\hline
       Symbol  & Quantity  & Symbol  &Quantity  \\\hline \hline
      $\rho$   & Density  &$T$   & Temperature \\\hline
      $\mu$ & Dynamic viscosity  & $\lambda$ & Thermal conductivity \\\hline
      $c$   & Adiabatic speed of sound & $c_T$ & non conductive speed of sound\\\hline
      $C_v$ & Isochoric heat capacity  &  $C_p$ & Isobaric heat capacity \\\hline
      $\gamma$ & Ratio  $C_p/C_v$  &  $\Gamma$ & Grüneisen parameter \\\hline
      $Pr$ & Prandtl number $\mu\,C_p/\lambda$   & $\text{Kn}$  & Knudsen number \eqref{Knudsen} \\\hline
    \end{tabular}
    \caption{Notation for a single fluid model (see also \eqref{CT} and \eqref{Cv_etc}).}
    \label{tab4}
\end{table}
In Tables \ref{tab2} and \ref{tab3}, the notation and relations used for two-fluid mixtures are given.

\begin{table}[ht!]
    \centering
    \begin{tabular}{|c|l||c|l|}\hline
       Symbol  & For fluid $\pm$ & Symbol  & For fluid $\pm$ \\\hline \hline
      $\alpha^\pm$   & Volume fraction  &$\eta^\pm$   & Mass fraction \\\hline
      $\rho^\pm$ & Density  & $v^\pm$ & Specific volume ($={1}/{\rho^\pm}$) \\\hline
      $c^\pm$   & Adiabatic speed of sound & $e^\pm$ & Specific internal energy\\\hline
      $s^\pm$ & Specific entropy  &  $E^\pm$ & Specific total energy \\\hline
      $p^\pm$ & Pressure  &  $H^\pm$ & Specific total enthalpy \\\hline
      $u^\pm$ & Velocity  &  $T^\pm$ & Temperature\\\hline
      $\mu^\pm$ & Viscosity  &  $h^\pm$  & Specific enthalpy \\\hline
       $\Gamma^\pm$ & Grüneisen parameter& $\text{Kn}^\pm$  & Knudsen numbers \eqref{Knudsen_pm} \\\hline
    \end{tabular}
    \caption{Notation for quantities per fluid.}
    \label{tab2}
\end{table}

\begin{table}[ht!]
    \centering
    \begin{tabular}{|l||l|}\hline
      $1=\alpha^++\alpha^-$   & $\eta\equiv\eta^+-\eta^- $   \\\hline
      $\eta^\pm=(1\pm\eta)/2\ge 0$   & $|\eta|\le 1=\eta^++\eta^- $   \\\hline
      $\rho=\alpha^+\,\rho^++\alpha^-\,\rho^-$ & $v=\eta^+\,v^++\eta^-\,v^-={1}/{\rho}$  \\\hline
      $e=\eta^+\,e^++\eta^-\,e^-$ & $s=\eta^+\,s^++\eta^-\,s^-$  \\\hline
      $E^\pm=e^\pm+\frac12|u^\pm|^2$   &$u=\alpha^+\,u^++\alpha^-\,u^-$   \\\hline
      $E=\eta^+\,E^++\eta^-\,E^-$   &$H=\eta^+\,H^++\eta^-\,H^-$   \\\hline
      $p=\alpha^+p^++\alpha^-p^-$ & $\mu=\alpha^+\,\mu^++\alpha^-\,\mu^-$  \\\hline
      $h^\pm=e^\pm+\frac{p^\pm}{\rho^\pm}$ & $H=E+\frac{p}{\rho}$  \\\hline
    \end{tabular}
    \caption{Notation for mean relations.}
    \label{tab3}
\end{table}

\subsection{On the speed of sound in single fluid}\label{sec_sos}

Following L. Landau and E. Lifchitz \cite{Landau}, a sound wave is a small amplitude oscillatory motion in a single compressible fluid. These authors derive the classical linear wave equation where the speed of sound is the usual thermodynamic coefficient which  
holds this name. Extending this definition to the mixture of two non miscible fluids leads us to consider the propagation of waves in such a medium. According to Whitham \cite{Whitham}, there is no single precise definition of what exactly constitutes a wave. Nevertheless this Author proposes to distinguish between hyperbolic waves, see Definition \ref{def1}, and dispersive waves. The latter being plane-waves where frequency $\omega$ is a defined real function of the wave number $k$ and the function $\omega (k)$ is determined by the particular system under consideration. In this case, the speed of the wave is its phase speed, that is $\omega (k)/k$, see (\ref{plane_wave}), and the waves are usually said to be "dispersive" if this phase speed is not a constant but varies with $k$.

\paragraph{Dispersion relation: the classical case of the Navier-Stokes equation}
Considering the classical compressible Navier-Stokes equation in $1D$, see \eqref{systeme_NST}, 
it is shown in Section \ref{AppenNS} that small disturbances around a constant state of rest will propagate as the superposition of plane waves ($k\in\mathbb{R}\,,\omega\in\mathbb{C}$):
\begin{equation}\label{plane_wave}
    W= W_0\,\exp{i(k\,x-\omega\,t)}=W_0\,\exp{ik\left(x-\frac{\omega}{k}\,t\right)}\,,
\end{equation}
provided $\omega$ and $k$ satisfy equation \eqref{NS1D} hereafter, which is known as the dispersion relation.
With the notation given in Table \ref{tab4}, setting:
\begin{equation}\label{CT}
a = \frac{2\,k\, \mu}{3\, \rho}\equiv k\,a_1\quad
\mbox{ and }\quad c_T^2  \equiv\frac{\partial p}{\partial \rho}\big|_T\,,
\end{equation}
where $a$ is Stokes' attenuation \cite{Stokes},
the dispersion relation \eqref{RD_NS} found in Section \ref{AppenNS} can be written as:
\begin{multline}
\label{NS1D}
\left(\frac{\omega}{k}\right)^{3} + 2\,i\,k\,a_1  \,\left( 1+\frac{3\,\gamma}{4\,Pr} \right)\left(\frac{\omega}{k}\right)^{2}  -  \left(c^{2} + \frac{3\,\gamma}{\text{Pr}}a_1^2\,k^2\right)\left(\frac{\omega}{k}\right) \\
- i\,\frac{3\,\gamma\,k\,a_1\,c_T^2}{2\,Pr} = 0\,,
\end{multline}

 In general this dispersion relation has $3$ solutions $(\omega\,,k)$ with $\omega\in\mathbb{C}$ and $k\in\mathbb{R}$\,. By writing $\omega=\omega(k)=\omega_R+i\,\omega_I$\,, it is seen that \eqref{plane_wave} reads:
\begin{equation}\label{plane_wave_bis}
    W= W_0\,(\exp{-\sigma(k)\,t)}\cdot\exp{ik\left(x-c(k)\,t\right)}\,,
\end{equation}
\begin{equation}\label{def_c(k)}
\mbox{ where }     c(k)\equiv \frac{\omega_R(k)}{k}\,,\mbox{ and } \sigma(k)\equiv -\omega_I(k)\,.
\end{equation}
Sticking to Whitham's definition of dispersive waves \cite{Whitham}, the phase speed of the wave is $c(k)$, while $\sigma(k)$ corresponds to the attenuation for positive values or amplification for negative values.

Although the dispersion relation \eqref{NS1D} could be explicitly solved for $\omega/k$ (with fixed frequency $k$) using Cardano's formula for third order polynomials, it is more intelligible to expand the $\omega/k$ roots with respect to $k$ for relatively small $k$, where the medium continuum hypothesis applies. More precisely, with respect to the Knudsen number based on $k$. See in Section \ref{AppenNS}, \eqref{Knudsen} and after.


\subsection{On the speed of sound in two-fluid mixtures: State of the art review}\label{sec142}

The goal of this work is to study the speed of sound in two-fluid non miscible mixtures. Following L. Landau and E. Lifchitz, small amplitude oscillatory motions in fluid mixtures are studied. Surprisingly, in a two-fluid mixture, the measured speed of sound can be one order of magnitude smaller than that of its constituents. For example for water and air in normal conditions the speed of sound in the mixture can be about $23\,\text{m/s}$ while it is $1\,500\,\text{m/s}$ in water and $330\,\text{m/s}$ in air, see Figure \ref{fig1j1}. This is supported by both analytical and experimental results, as will now be discussed.
\subsubsection{Analytical formulas from the literature}
To determine the speed of sound in a mixture of two-fluids, the most common and basic model is based on the mixture effective mean compressibility and mean density (Wood \cite{wood}, Brennen \cite{Brennen}). As in the literature, {\it Wood's speed of sound} will refer to the sound speed obtained by this model, and denote it by ${\bm c_w}$. It can be straightforwardly derived as follows. Consider a homogeneous mixture of two non miscible fluids in mechanical equilibrium (equal pressure), and assume that the mass fraction of each fluid is constant. Then differentiating the relation ${1}/{\rho}=\eta^-/\rho^-+\eta^+/\rho^+$
with respect to the pressure leads to:
\begin{equation}\label{wood0}
\frac{1}{\rho^2 c^2_w} = \frac{\eta^-}{(\rho^-) ^2 (c^-) ^2}+\frac{\eta^+}{(\rho^+)^2  (c^+) ^2}\,.
\end{equation}
This gives Wood's formula and is also equivalent to (see Figure \ref{fig1j} for a plot of $c_w$ against $\alpha^+$):
\begin{equation}\label{wood1}
    \frac{1}{\rho c^2_w} = \frac{\alpha^-}{\rho^- (c^-) ^2}+\frac{\alpha^+}{\rho^+  (c^+) ^2}\,.
\end{equation}

As will be shown later in Section \ref{sec311}, this derivation is equivalent to considering a homogeneous medium where the sound propagates with no slipping of thermal exchange between the two phases, so out of thermal equilibrium between the two fluids (Stewart and Wendroff \cite{Stewart}).

Other models for the speed of sound in two-fluid mixtures are presented in Nguyen {\it et al.} \cite{nguyen}. These models depend on the flow regime under consideration (homogeneous, stratified and slug regimes), and use the consideration that the interface of one phase is a confining elastic boundary for the other, making use of the general theory of sound propagation in fluids confined by elastic tubes.
\begin{itemize}
    \item[(i)] For a stratified flow regime, the model assumes the existence of two sound speeds, one for each fluid ($c_{eff}^{\pm}$), and gives the following formulas for the effective speed of sound  in each phase as a function of the pure fluid velocities $c^{\pm}\,$:
\begin{equation}\label{nguyen1}
\frac{1}{(c_{eff}^{-})^2} = \frac{1}{(c^{-})^2} + \frac{\alpha^+}{\alpha^-}\frac{\rho^-}{\rho^+} \frac{1}{(c^{+})^2}\,,
\end{equation}
\begin{equation}\label{nguyen2}
\frac{1}{(c_{eff}^{+})^2} = \frac{1}{(c^{+})^2} + \frac{\alpha^-}{\alpha^+}\frac{\rho^+}{\rho^-} \frac{1}{(c^{-})^2}\,.
\end{equation}

In each of the two previous formulas the second term is due to the confinement of medium $\mp$ by the other elastic medium $\pm$.
\item[(ii)] For a slug flow regime, and assuming that the sound wave passes through each phase successively, the model proposed in Nguyen {\it et al.} \cite{nguyen} gives the following formula for the speed of sound :
\begin{equation}\label{nguyen3}
\frac{1}{c_{eff}^{+}} = \frac{\alpha^+}{c^{+}} + \frac{\alpha^-}{c^{-}}\,.
\end{equation}
\item[(iii)] For a homogeneous bubbly flow regime, the model proposed in Nguyen {\it et al.} \cite{nguyen} takes into account the elastic confinement by the other phase and the following two formulas for the speed of sound in each phase are proposed :
\begin{equation}\label{cgeff}
\frac{1}{(c_{eff}^{-})^2} = \frac{\alpha^-}{(c^{-})^2} + \frac{\rho^-}{\rho^+} \frac{\alpha^+}{(c^{+})^2}\,,
\end{equation}
and 
\begin{equation}\label{cweff}
\frac{1}{(c_{eff}^{+})^2} = \frac{\alpha^+}{(c^{+})^2} + \frac{\rho^+}{\rho^-} \frac{\alpha^-}{(c^{-})^2}\,.
\end{equation}
The composite speed of sound in the medium is derived assuming that the sound wave passes through each phase successively :
\begin{equation}
\label{nguyen-bulles}
\frac{1}{c} = \frac{\alpha^+}{c_{eff}^{+}}  + \frac{\alpha^-}{c_{eff}^{-}}\,.
\end{equation}
\end{itemize}

\begin{figure}[ht!]
\begin{center}
\includegraphics[scale=0.6]{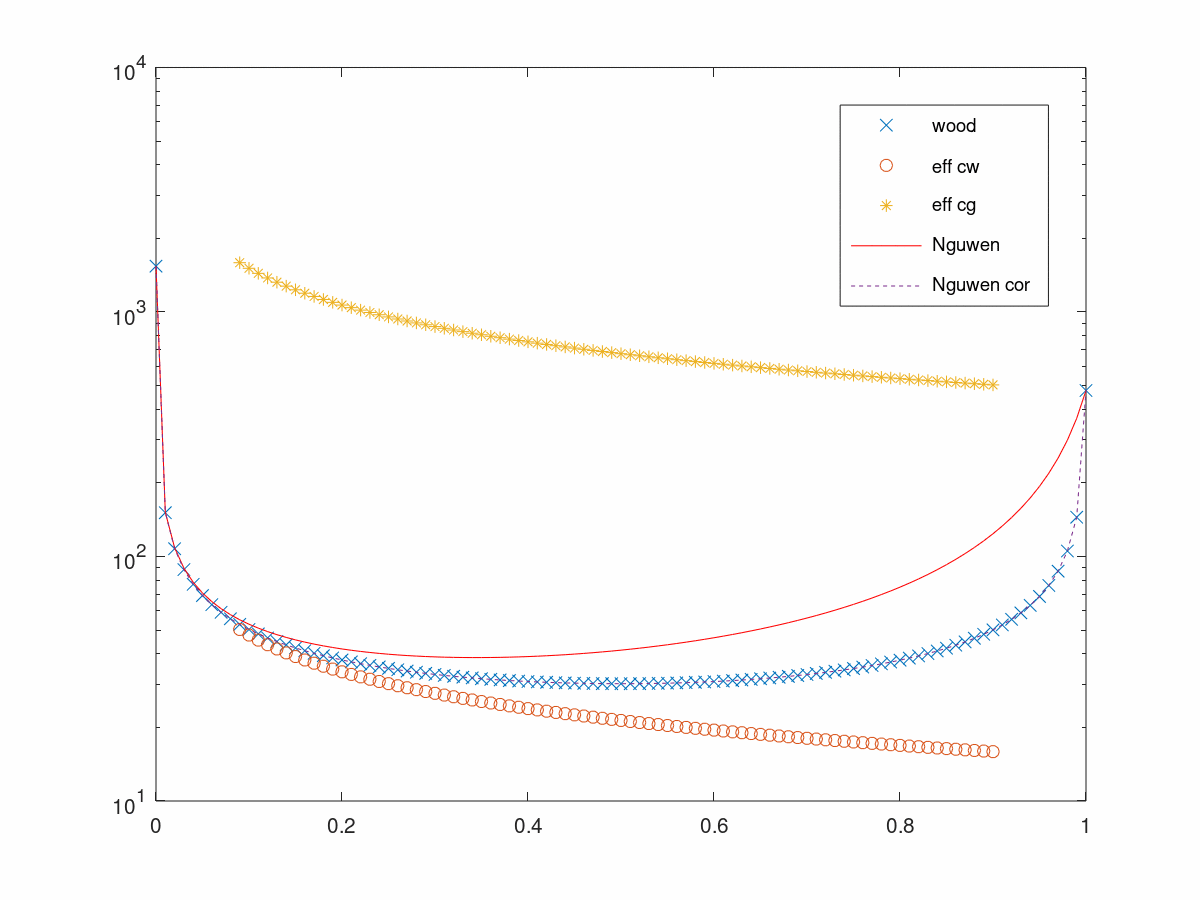}
\caption{Comparison of different sound speed relations for a water-air bubbly mixture as a function of void fraction $\alpha^-$.}
\label{fig23}
\end{center}
\end{figure}

The hypothesis that the sound wave passes through each phase successively in a bubbly homogeneous flow is not necessary applicable and experimental work shows a poor agreement with (\ref{nguyen-bulles}) (see section \ref{sec:exp}). A more relevant approach is to consider the mixture compressiblity based on the two effective compressibilities of the fluids computed by (\ref{cgeff}) and (\ref{cweff}):
\begin{equation}
    \frac{1}{\rho c^2} = \frac{\alpha^-}{\rho^- (c_{eff}^-) ^2}+\frac{\alpha^+}{\rho^+  (c_{eff}^+) ^2}\, \\
    = \frac{\alpha^-}{\rho^- (c^-) ^2}+\frac{\alpha^+}{\rho^+  (c^+) ^2}\,.
\end{equation}
Note that this is Wood's formula (\ref{wood1}). This means that the funneling effects, of each phase on the sound wave propagation in the other phase, cancel each other out, for the mean sound speed in the medium.
Figure \ref{fig23} shows a plot of the different formulas above for water and air bubbly mixtures as a function of void fraction.

Generalizations of Wood's formula that account for different additional phenomenon, such as surface tension, phase change in a single component mixture (e.g. water and its vapor), can be found in the books of Brennen \cite{Brennen} and Landau and Lifshitz \cite{Landau}. 

Temkin in \cite{Temkin92} refers to the previous Wood's formula as the frozen equilibrium speed of sound, because it is derived without taking into account thermal exchanges between the two phases. Temkin gives a relation for a complete equilibrium speed of sound that is equivalent to the formula put forward here in \ref{thermo:equilibrium}.

The models cited above do not consider the effect of the wave frequency and hence the sound dispersion and attenuation. The first model that considered the effect of frequency on sound speed propagation in bubbly liquids is due to Foldy \cite{Foldy} and is also presented in Feshbach {\it et al.} \cite{Feshbach2} page 1498. This model is based on a multiple scattering analysis. 
If $+$ is the liquid phase and $-$ is the bubbly gas phase, the speed of sound in a bubbly medium, with bubbles of radius $R$ is given by:
\begin{equation}\label{foldy0}
\frac{1}{c^2} = \frac{k^2(\omega)}{\omega^2} = \alpha^+\,\rho^+ \left( \frac{1}{\rho^+ (c^+)^2} + \frac{\alpha^-}{\rho^- (c^-)^2} \frac{1}{1 - \frac{R^2 \rho^+ \omega^2}{3 \rho^- (c^-)^2} (1+iR\frac{\omega}{c^+})} \right)\,,
\end{equation}
where here $c=\omega/k$ is a complex number, see \eqref{plane_wave_bis}-\eqref{def_c(k)}.

Noting that $\omega_0= \frac{c^-}{R}\sqrt{\frac{3\, \rho^-}{\rho^+}} $ is the Minnaert resonance frequency, we have:
\begin{equation}\label{foldy}
\frac{1}{c^2} = \frac{k^2(\omega)}{\omega^2} = \alpha^-\,\rho^+ \left( \frac{\alpha^-}{\rho^+ (c^+)^2} + \frac{\alpha^+}{\rho^- (c^-)^2} \frac{1 }{1 - \frac{\omega^2}{3 \omega_0^2} (1+i R \frac{\omega}{c^+})} \right).
\end{equation}

\begin{figure}[ht!]
\begin{center}
\includegraphics[scale=0.3]{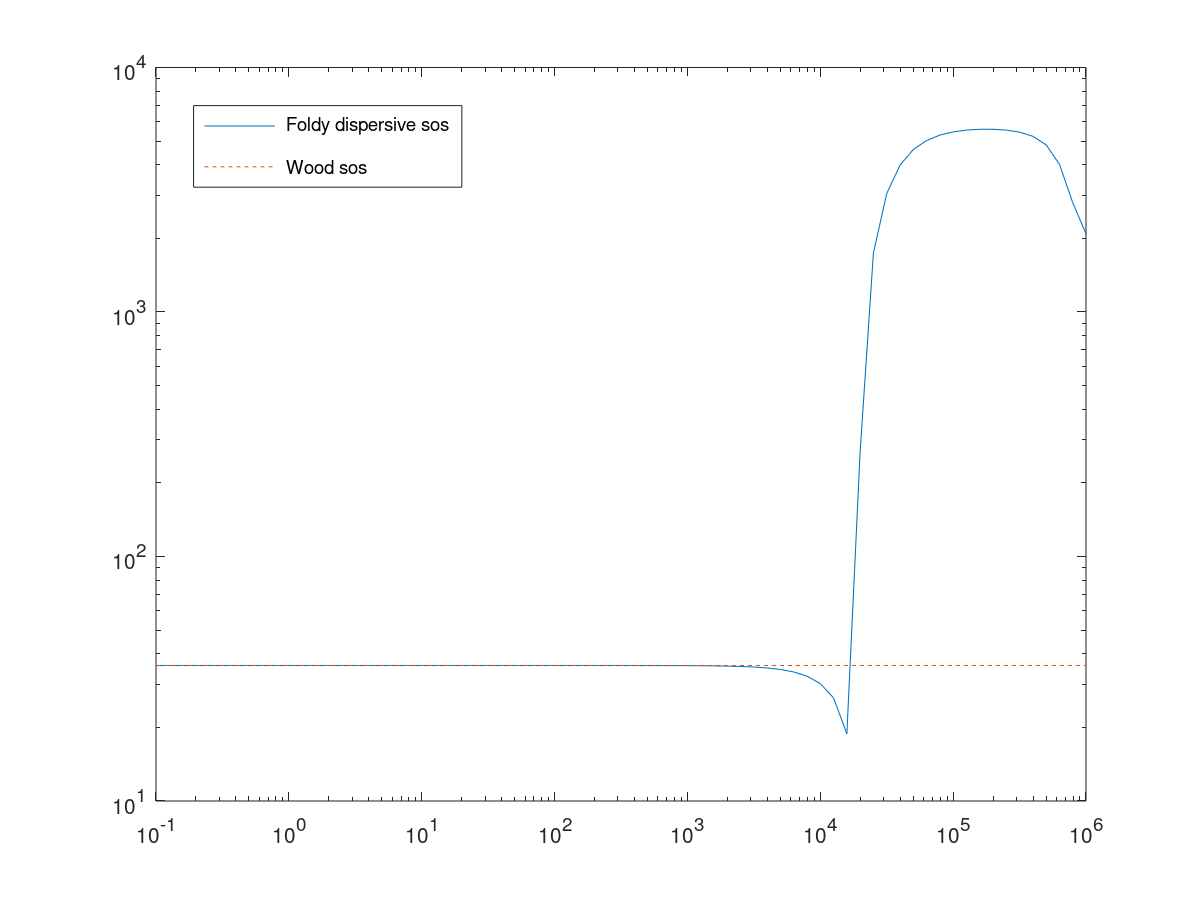}
\caption{Speed of sound (sos) in equation (\ref{foldy}), for a water-air mixture at a void fraction of $\alpha^-= 0.1$, and for a range of wave frequencies.} \label{sos_water_air_foldy}
\end{center}
\end{figure}

This model has been generalized to single substance two-fluid mixtures ({\it i.e.} with phase change) by Trammell \cite{Trammell} and Kielland \cite{Kielland}:
\begin{equation}\label{Trammell}
\frac{k^2(\omega)}{\omega^2} = \frac{\alpha^-}{1+2\alpha^+}\rho^+ \left( \frac{1}{\rho^+ (c^+)^2} + \frac{\alpha^+}{K_b} \frac{1 }{1 - \frac{R^2 \rho^+ \omega^2 K_b}{3} (1+ik(\omega)R)}\right),
\end{equation}
with $K_b = \frac{1}{\rho^- (c^-)^2} + \frac{3\lambda^+}{\rho^- L \frac{dP}{dT}} \left[ \frac{1}{R}\sqrt{\frac{C_p^+}{\lambda^+ \omega}} + \frac{i}{R^2\omega}\right]$ and where $\lambda^+$ is the liquid thermal conductivity, $C_p^+$ its isobaric heat capacity, $L$ the latent heat of evaporation and $\frac{dP}{dT}$ is the pressure derivative along the saturation curve. 
Figure \ref{sos_water_air_foldy} shows the speed of sound from equation (\ref{foldy}), for a water-air mixture at a void fraction of $\alpha^-= 0.1$, and for a range of wave frequencies.

Another theoretical model for the sound dispersion and attenuation is due to Temkin \cite{Temkin00} and this gives predictions numerically close to Foldy's model \cite{Temkin00}.

Substantial experimental work has been conducted to validate the previously mentioned formulas.

\subsubsection{Experiments at low frequency}\label{sec:exp}
Many authors have produced experiments studying the variation of the speed of sound with respect to the void fraction in either a two-substance two-phase system (e.g. air and water) or in a single-substance two-phase system (e.g. water and steam). Among the works that do not take into account the effect of the frequency, the reader is referred to Costigan {\it et al.} \cite{Costigan97} and Hiva {\it et al.} \cite{Hiva10}. Their measurements, as well as those of Karplus  \cite{Karplus61}
show a good agreement with Wood's formula, but with the additional hypothesis that the gas obeys a perfect-gas equation of state with a polytropic index $\gamma =1$, corresponding to an non-heat-conducting speed of sound $c^-_T$ in the gas:
\begin{equation}\label{c_woodT}
    \frac{1}{\rho \tilde{c}_w^2} = \frac{\alpha^-}{\rho^- (c^-_T) ^2}+\frac{\alpha^+}{\rho^+  (c^+) ^2}\,.
\end{equation}

\begin{figure}[ht!]
\label{Karplus}
\begin{center}
\includegraphics[scale=0.7]{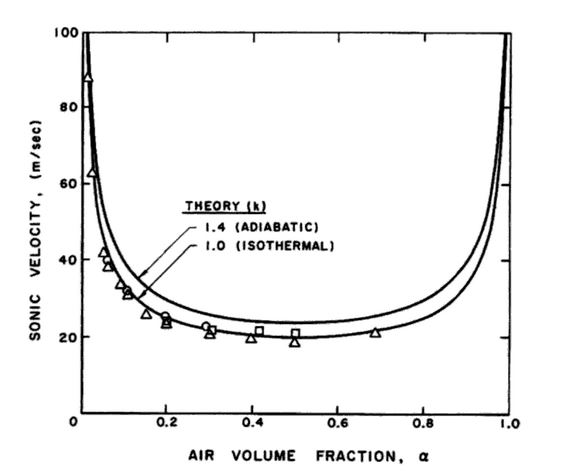}
\caption{Karplus \cite{Karplus61} data, taken from Brennen \cite{Brennen}. The isentropic Wood formula overestimates the experimental speed of sound. Using a non conductive sound speed for air gives a better match with the experimental results.}
\label{figkarplus}
\end{center}
\end{figure}

Costigan {\it et al.} \cite{Costigan97} also show that Nguyen's formula for bubbly liquids (\ref{nguyen-bulles}) gives a less precise prediction. 


Similarly, Temkin \cite{Temkin00} shows that the (isentropic) Wood's formula overestimates slightly the speed of sound (as can be seen in Figure \ref{figkarplus}) and instead of modifying the polytropic index for the gas, it is argued that the relaxed equilibrium speed of sound (see section \ref{thermo:equilibrium}) better matches the data from Karplus \cite{Karplus61} and other experiments, and hence corresponds to the best choice of model to adopt.

\subsubsection{Experiments at varying frequencies}

\begin{figure}[ht!]
\begin{center}
\includegraphics[scale=0.7]{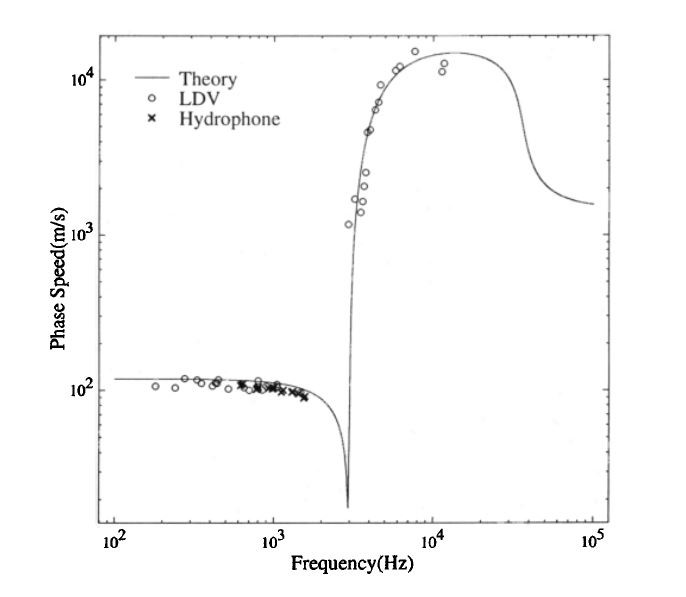}
\caption{Speed of sound versus frequency from \cite{Cheyne95}. Void fraction $= 1\%$ and bubble radius $= 1.11mm$. Taken from \cite{Cheyne95}.}
\label{Cheyne-exp}
\end{center}
\end{figure}

Little experimental data are available for the dependence of speed of sound (and sound attenuation) on the wave frequency. The available data in the literature, in Silberman \cite{Silberman57} and in Cheyne {\it et al.} \cite{Cheyne95}, show a good agreement with the scattering theory, as shown in Figure \ref{Cheyne-exp}.

Concerning the use of numerical simulation to evaluate the speed of sound in a mixture of two non miscible fluids, the reader is referred to the papers by J.A. Redford {\it et al} \cite{Redford} and V. Leroy {\it et al} \cite{Leroy2}.



\subsection{A few plots for the Speed of Sound}
In the present article, analytical formulas for the speed of sound will be derived based on various dynamical models. 
These formulas are gathered in Table \ref{tab1}. One of the major outcomes of this work concerns the three basic formulas: (i) $c_0$ found in Section \ref{sec521}, (ii) $c_\kappa$ found in Section \ref{sec323}  and (iii) $c_w$ found in Section \ref{sec311}. In the first two cases the two-fluids can slip freely (two independent velocities), while in the third case the fluids have the same velocity.\\
In general, when the velocities are independent there is at least a drag force that tends to equalize these velocities. The speed of sound is unchanged so long as this restoring force does not involve derivatives of the dependent variables, which is in contrast to the added mass force that will now be discussed.\\
The values $c_0$ and $c_w$ correspond to the extremes ({\it i.e.} asymptotic cases) where the physical reality is somewhere between. Introducing an added mass force in the evolution equations for the velocities (momentum equations), see \eqref{f1}, it is indeed possible to interpolate between these two extreme cases and this leads to an intermediate speed of sound $c_\kappa$. Here $\kappa$ is a non-dimensional coefficient that governs the strength of this restoring force whose effect tends to equalize the velocities. Actually, as expected and shown by \eqref{interpol}, $c_0$ corresponds to $c_\kappa$ for the value $\kappa=0$ ({\it i.e.} no restoring force) and $c_w$ corresponds to the limit of $c_\kappa$ as $\kappa$ tends to infinity (equal velocities).\\

\begin{figure}[ht!]
\begin{center}
\includegraphics[scale=0.3]{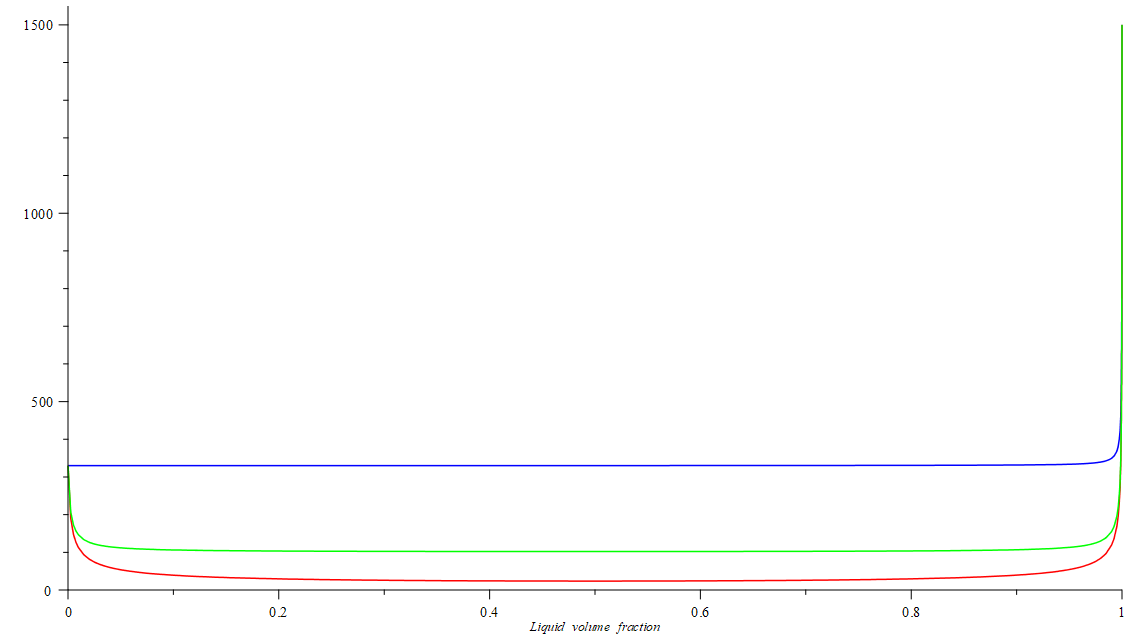}
\caption{Sound speeds, $c_w$ (red), $c_\kappa$ (green), $ \kappa=10$,  and $c_0$ (blue), as a function of water volume fraction $\alpha^+$. Enlarged plots for small or large $\alpha^+$ are given Figure \ref{figzoom} and a plot with an adapted $\alpha^+$ scale to enhance readability is given in Figure \ref{fig1j1}.} \label{fig1j}
\end{center}
\end{figure}

For the air and water mixture shown in Figures \ref{fig1j} and \ref{fig1j1}, it is seen that $c_0$ is a monotone function of the water liquid volumic fraction $\alpha^+$\,, while $c_w$ is decreasing to a minimum and then increases. The same behavior is observed for $c_\kappa$ (for sufficiently large $\kappa$). This matches with the experimental observation mentioned in Section \ref{sec142}.
\begin{table}[ht!] 
    \centering
    \begin{tabular}{|c|l|l|}\hline
       Notation  & Denomination & Value \\\hline \hline
      $c$   & Thermodynamic isentropic SoS & (\ref{eq:thermo_c2_s}) \\\hline
      $c_T$ & Thermodynamic non conductive SoS & (\ref{eq:thermo_c2_T}) \\\hline
      $c_w$ & Wood formula &  (\ref{wood1}) \\\hline
      $c_\kappa$ & Added mass model SoS &  (\ref{SoS1}) \\\hline
      $c_0$ & Two-velocity model SoS &  (\ref{vitson2vit}), (\ref{SoS2}) \\
      \hline
     $c_{eq}$ & Equilibrium/one temperature SoS & (\ref{cS-equilibrium}) \\
      \hline
    \end{tabular}
    \caption{Expressions for the various Speed of Sound (SoS) in this article.}
    \label{tab1}
\end{table}

\begin{figure}[ht!]
  \begin{minipage}[b]{0.45\linewidth}
   \centering
   \includegraphics[scale=0.3]{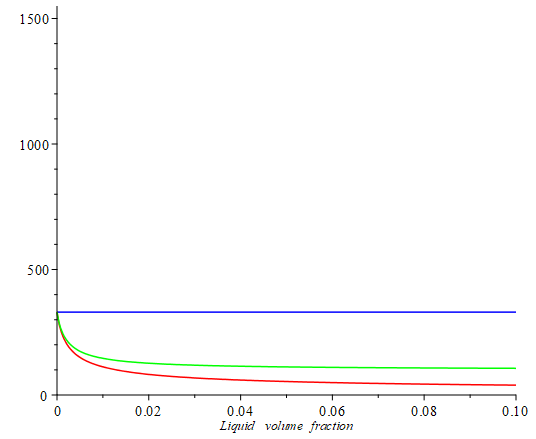}    
  \end{minipage}
\hfill
  \begin{minipage}[b]{0.55\linewidth}
   \centering
   \includegraphics[scale=0.3]{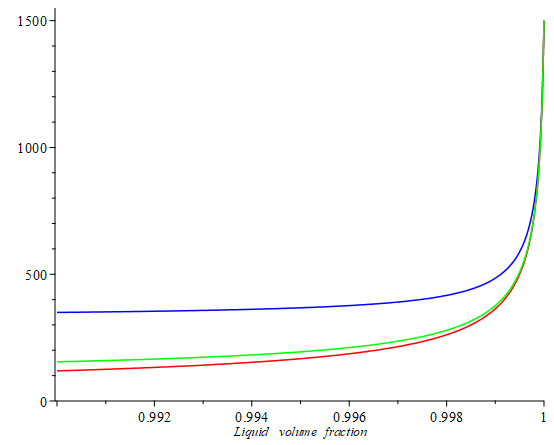}     
  \end{minipage}

  \caption{Behavior of $c_w$ (red), $c_\kappa$ (green), $ \kappa=10$,  and $c_0$ (blue) for small or large water volume fraction $\alpha^+$.}
  \label{figzoom}
\end{figure}
\begin{figure}[ht!]
\begin{center}
\includegraphics[scale=0.3]{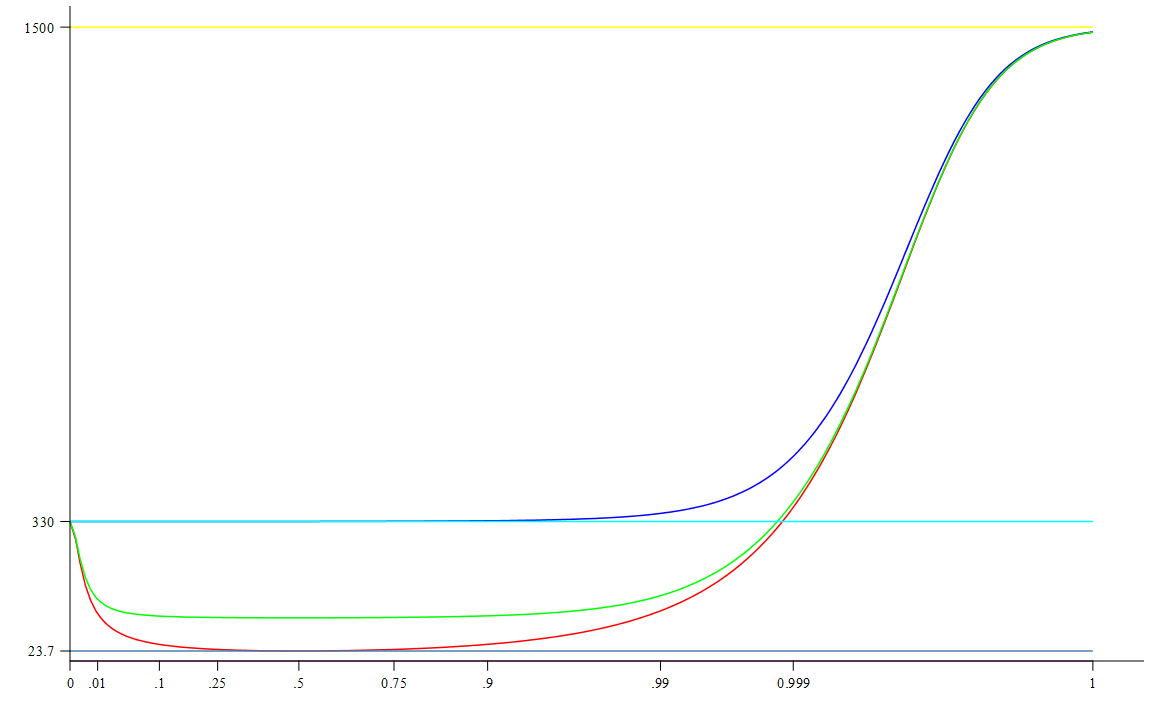}
\caption{Sound speeds, $c_w$ (red), $c_\kappa$ (green), $ \kappa=10$,  and $c_0$ (blue), as a function of the water volume fraction $\alpha^+$. The scale of the $x$-axis is adapted to the variations of these functions in the vicinity of $\alpha^+=1$ (pure water, yellow) and $\alpha^+=0$ (pure air, turquoise). The level corresponding to the minimal speed of sound is also plotted (marine blue). The parametrization $\theta$ of the $x$-axis is $0\le \theta < +\infty$ with $\alpha^+=1-{10}^{-{\theta}^{2}}$. } \label{fig1j1}
\end{center}
\end{figure}

\section{two-fluid flow models}\label{sec_2-fluid_mod}
In this Section, the main models encountered in the literature are briefly reviewed. Roughly speaking, there are two classes of model. In the first, local disequilibrium between the two fluids is ignored, while in the second this hypothesis is not made. Recalling that these models are averaged models, hence, average variables at a given point $(x\,,t)$ are considered. Let us take, for instance, the relative velocity $u^+-u^-$. Assuming that $u^+=u^-$ at a point means that on average there is no slipping between the two fluids (or that it is neglected in the model). With respect to a two-velocity model, this is an asymptotic model that can be relevant, or not, depending on the physical situation. In this Section, inviscid (or perfect) fluids will not be considered. The effect of viscosity and thermal diffusivity will be addressed in Section \ref{Sec4}.\\

From the mathematical point of view, well-posedness is, up to now, an outstanding open question and the reader is referred to Bresch {\it et al.} \cite{Bresch} concerning its discussion. The two difficulties are (i) the presence of non conservative terms {\it e.g.} the term $ p\,\alpha^\pm_t$ in the energy equation (\ref{1u3}) or $\alpha^\pm\nabla p$ in the momentum equation \eqref{1u1T2} and (ii) non-hyperbolicity of the convective part of the equation (Stewart and Wendroff \cite{Stewart}). The question of well-posedness will not be addressed in this paper as the focus will be on determination of the speed of sound, which deals with small perturbations around a rest state with vanishing velocities.\\

Listed below are the different inviscid models considered in this paper, referring to Section \ref{Sec3} for the expression of the associated speed of sound. The definitions of the quantities appearing in the different systems are given in Tables~\ref{tab2}~and~\ref{tab3}.

\subsection{One-velocity models}
In these models, it is assumed that the two fluids are in mechanical equilibrium $(p = P^{\pm}(\rho^{\pm},s^{\pm}))$ and that there is no slipping between the fluids $u=u^{\pm}$.
\subsubsection{Isentropic case}
The two-fluid one-velocity system of equations reads:
\begin{equation}\label{isoS1u1}
(\alpha^\pm\rho^\pm)_t+\text{div}(\alpha^\pm\rho^\pm u)=0\,,
\end{equation}
\begin{equation}\label{isoS1u2}
(\rho u)_t+\text{div}(\rho u\otimes u)+ \nabla p=0\,,
\end{equation}
where $a\otimes b$ denotes the matrix $(a\otimes b)_{i\,j}=a_i\,b_j$ and the two EoS are:
\begin{equation}\label{isoS1u3}
    EoS_\pm(p,\rho^\pm)=0\,.
\end{equation}

\subsubsection{A model with two temperatures}\label{sec:1vit}

Here it is again assumed that the two-fluids are in mechanical equilibrium and that there is no slipping between the two fluids, $u = u^{\pm}$, but two energy equations are considered in the model:
\begin{equation}\label{1u1}
(\alpha^\pm\rho^\pm)_t+\text{div}(\alpha^\pm\rho^\pm u)=0\,,
\end{equation}
\begin{equation}\label{1u2}
(\rho u)_t+\text{div}(\rho u\otimes
u)+ \nabla p=0\,,
\end{equation}
\begin{equation}\label{1u3}
(\alpha^\pm\rho^\pm E^\pm)_t+\text{div}(\alpha^\pm(\rho^\pm E^\pm+p) u)+
 p\alpha^\pm_t=0\,,
\end{equation}
\begin{equation}\label{1u4}
p = P^{\pm}(\rho^{\pm},s^{\pm}).
\end{equation}

For smooth solutions this system is equivalent to 
\begin{equation}\label{1u1_bis}
\rho_t+div(\rho\, u)=0\,,
\end{equation}
\begin{equation}\label{1u2_bis}
(\rho u)_t+\text{div}(\rho u\otimes u)+ \nabla p=0\,,
\end{equation}
\begin{equation}\label{1u3_bis}
\eta_t+ u.\nabla\eta =0\,,
\end{equation}
\begin{equation}\label{1u4_bis}
s^{\pm}_t+ u.\nabla s^{\pm} =0\,,
\end{equation}
\begin{equation}\label{1u5_bis}
    p = P(\rho,\eta,s^+,s^-)\,.
\end{equation}


\subsubsection{A model with thermal equilibrium}

It is supposed that in addition to mechanical equilibrium between the two fluids ($p^+ = p^-$), the thermal equilibrium $T^+ = T^-$ is also verified at all times. The corresponding dynamic model is as follows:
\begin{equation}\label{1u1T1}
(\alpha^\pm\rho^\pm)_t+\text{div}(\alpha^\pm\rho^\pm u)=0\,,
\end{equation}
\begin{equation}\label{1u1T2}
(\alpha^\pm\rho^\pm u)_t+\text{div}(\alpha^\pm\rho^\pm u\otimes
u)+ \alpha^\pm\nabla p=0\,,
\end{equation}
\begin{equation}\label{1u1T3}
(\rho E)_t+\text{div}((\rho E +p) u) =0\,,
\end{equation}
\begin{equation}\label{1u1T4}
p = p^{\pm}(\rho^{\pm},e^{\pm}) \,,\,  T = T^{\pm} (\rho^{\pm},e^{\pm})\,.
\end{equation}


For smooth solutions, the above system is equivalent to   

\begin{equation}\label{1u1T1_bis}
\rho_t+\text{div}(\rho u)=0\,,
\end{equation}
\begin{equation}\label{1u1T2_bis}
(\rho u)_t+\text{div}(\rho u\otimes u)+ \nabla p=0\,,
\end{equation}
\begin{equation}\label{1u1T3_bis}
\eta_t+ u.\nabla \eta =0\,,
\end{equation}
\begin{equation}\label{1u1T4_bis}
s_t+ u.\nabla s =0\,,
\end{equation}
\begin{equation}\label{1u1T5_bis}
    p = P(\rho,\eta,s)\,.
\end{equation}

\subsection{Two-velocity models}

\subsubsection{Isentropic case}
In the isentropic case the two-velocity model reads as follows:

\begin{equation}\label{2u1}
(\alpha^\pm\rho^\pm)_t+\text{div}(\alpha^\pm\rho^\pm u^\pm)=0\,,
\end{equation}
\begin{equation}\label{2u2}
(\alpha^\pm\rho^\pm u^\pm)_t+\text{div}(\alpha^\pm\rho^\pm u^\pm\otimes
u^\pm)+ \alpha^\pm\nabla p=0\,,
\end{equation}

 where
 the pressure $p$ is related to the densities $\rho^\pm$
through two equations of state
\begin{equation}\label{2u3}
    EoS_\pm(p,\rho^\pm)=0\,.
\end{equation}

\subsubsection{A model with two temperatures}

The system (\ref{2u1})-(\ref{2u2}) can be completed with the two energy equations:
\begin{equation}
\label{2u4}
\frac{\partial (\alpha^\pm \rho^\pm E^\pm)}{\partial t} + \text{div}\,\, (\alpha^\pm (\rho^\pm E^\pm+p) u^\pm)
+ p \frac{\partial \alpha^\pm}{\partial t} = 0,
\end{equation}
and (\ref{2u3}) should then be replaced with energy dependent equations of state.
\begin{equation}\label{2u3_bis}
EoS_\pm(p,\rho^\pm, s^\pm)=0\,.
\end{equation}

\subsubsection{A model with two temperatures and added mass}
The added mass modeling term is also called the ``virtual mass". Here, the approach found in the Appendix of Redford {\it et al.} \cite{Redford} is followed.\\
These models can be written as~:
\begin{equation}
\label{1.8}
\frac{\partial (\alpha^\pm \rho^\pm)}{\partial t} +
 \text{div}\,\, (\alpha^\pm \rho^\pm u^\pm) = 0,
\end{equation}
\begin{equation}
\label{1.10}
\frac{\partial (\alpha^\pm \rho^\pm u^\pm)}{\partial t} +
\text{div}\,\,\, (\alpha^\pm (\rho^\pm u^\pm \otimes u^\pm +  p{\bf I}))
= p\nabla \alpha^\pm+f^\pm\,,
\end{equation}
\begin{equation}
\label{1.12}
\frac{\partial (\alpha^\pm \rho^\pm E^\pm)}{\partial t} +
\text{div}\,\, (\alpha^\pm (\rho^\pm E^\pm+p) u^\pm)
+ p \frac{\partial \alpha^\pm}{\partial t} = f^\pm\cdot u^\pm\,,
\end{equation}

where
\begin{equation}
\label{f1}
f^\pm\equiv\mp\kappa\frac{\alpha^+\alpha^-\rho^+\rho^-}{\alpha^+\rho^++\alpha^-\rho^-}\frac{\partial (u^+-u^-)}{\partial t}\,,
\end{equation}
are the so-called added mass terms see {\it e.g.} M. Ishii and T. Hibiki  \cite{Ishii}.
The ten (in $3D$) differential equations \eqref{1.8} to \eqref{1.12} are supplemented with the two equations of state~:
\begin{equation}
\label{}
EoS_\pm(p,\rho^\pm, s^\pm)=0\,,
\end{equation}
and $\kappa\,(\ge0)$ may depend on thermodynamic variables ({\it i.e.} $p\,,s^\pm\,$) and $\alpha^\pm$.

\section{Kinematics and Speed of Sound}\label{Sec3}
In this Section the expressions for the speed of sound for each of the models presented in the previous Section are derived. As presented in the Introduction, the technique is to obtain the dispersion relation through plane-waves, as will now be discussed.\\

In $1D$, when linearized around a constant state, all the differential models given in Section \ref{sec_2-fluid_mod} (and also the Baer-Nunziato Model in Section \ref{BNModel}) lead to a linear differential system of the form:
\begin{equation}
\label{linearsystem}
\frac{\partial W}{\partial t} + A\, \frac{\partial W}{\partial x} = 0\,,
\end{equation}
where $A$ is a $N\times N$ constant matrix and $N$ is the number of differential equations occurring in the model.
\begin{Def}\label{def1}
The model \eqref{linearsystem} is said to be hyperbolic if there exists a basis $(r_1\,,\ldots\,,r_N)$ of $\mathbb{R}^N$ made with eigenvectors of $A$\,:
\begin{equation}
A\, r_k = \lambda_k\,r_k\,,\quad k=1\,,\ldots, N\,,
\end{equation}
where $\lambda_k\in \mathbb{R}$ are the associated eigenvalues.
\end{Def}
Dimensional analysis shows that the dimension in S.I. units of $\lambda_k$ is $\text{m/s}$,
{\it i.e.} velocity. Looking for non-vanishing plane-wave solutions, (\ref{plane_wave}), leads to the simple dispersion relation:
\begin{equation}\label{plane_wave1}
    \frac{\omega}{k}\in \Lambda\equiv\{\lambda_1\,,\ldots\,,\lambda_N\}\,,
\end{equation}
{\it i.e.} all these models are non dispersive equations since ${\omega}/{k}$ is constant.\\

When dealing with all the models in this paper,
and their linearization around a quiescent state,
it appears, for symmetry reasons, that the set $\Lambda$ is symmetric around $0$\,: $\Lambda=-\Lambda$\,. Hence $\Lambda=\{-\chi_\ell\,,\chi_\ell\,,\ell=1\,,\ldots, M\}$\,, $\chi_\ell\ge 0$. In general (but not always, see (\ref{2u1})-(\ref{2u2})) one of the $\chi_\ell$ vanishes. The velocity $0$ is in general associated with a material wave, such as the propagation of a contact discontinuity. Non-vanishing velocities $\chi_\ell$ correspond to pressure-waves and, besides the Baer-Nunziato Model (\ref{BN1})-(\ref{BN7}), see (\ref{LambdaBM}), it is found that for all the other models either $\Lambda=\{-\chi\,,\chi\,\}$ or $\Lambda=\{-\chi\,,0\,,\chi\,\}\,,\chi> 0$.

\subsection{One-velocity models}

\subsubsection{Isentropic case}\label{sec311}

In order to linearize the system (\ref{isoS1u1})-(\ref{isoS1u3}), it is noted that this is equivalent to:
\begin{equation}\label{isoS1u1_bis}
\rho_t+div(\rho u)=0\,,
\end{equation}
\begin{equation}\label{isoS1u2_bis}
(\rho u)_t+div(\rho u\otimes u)+ \nabla p=0\,,
\end{equation}
\begin{equation}\label{isoS1u3_bis}
\eta_t+ u.\nabla \eta =0\,,
\end{equation}
\begin{equation}\label{isoS1u4_bis}
    p = P(\rho,\eta)\,,
\end{equation}
for smooth solutions. Hence, a single-fluid Euler system is recognized, where the mass fraction variable $\eta$ takes the place of the usual entropy variable $s$.

The system (\ref{isoS1u1_bis}-\ref{isoS1u4_bis}) can be linearized around a constant $(\rho_0\,, u_0 \,, \eta_0)$ solution to obtain a system of type (\ref{linearsystem}) with $W = ~{}^t( \rho \,, u  \,, \eta )$ and 
\begin{equation}\label{A_1}
A = 
\left(
\begin{array}{ccc} 
u_0 & \rho_0 & 0  \\
\frac{1}{\rho_0} \left.\frac{\partial p}{\partial \rho }\right|_{\eta}& u_0 & \frac{1}{\rho_0}\left.\frac{\partial p}{\partial \eta }\right|_{\rho} \\
0 & 0 & u_0 
\end{array}
\right).
\end{equation}

The eigenvalues for the linearized system are $u_0$, $u_0-c$, $u_0$, where (dropping the subscripts $0$) the speed of sound $c$ is given by:
\begin{equation}\label{c2uneVit}
c^2=\frac {\rho^+\rho^-{(c^+)}^{2}{(c^-)}^{2}}{(\alpha^+\rho^++\alpha^-\rho^-)\left(\alpha^+\rho^-(c^-)^{2}+\alpha^-\rho^+(c^+)^{2}\right)}\,.
\end{equation}
which is equivalent to Wood's relation (using $\rho=\alpha^+\rho^++\alpha^-\rho^- \,, \eta^+ = \frac{\alpha^+ \rho^+}{\rho} $):
\begin{equation}\label{c_x}
\frac{1}{\rho^2 c^2}=\frac{\eta^+}{(\rho^+)^2 (c^+)^2} +\frac{\eta^-}{(\rho^-)^2 (c^-)^2} \, \mbox{ or } \frac{1}{\rho c^2}=\frac{\alpha^+}{\rho^+ (c^+)^2} +\frac{\alpha^-}{ \rho^- (c^-)^2}.
\end{equation}
The minimal speed of sound is attained for:
\begin{equation}\label{argmin_alphauneVit}
\alpha^+=\alpha_{min}=\frac{1}{2}\,{\frac {(\rho^+)^{2}{(c^+)}^{2}+(\rho^-)^{2}{(c^-)}^{
2}-2\,\rho^-\rho^+{(c^+)}^{2}}{-\rho^+\rho^-{(c^-)
}^{2}+(\rho^-)^{2}{(c^-)}^{2}+(\rho^+)^{2}{(c^+)}^{2}-\rho
^-\rho^+{(c^+)}^{2}}}\,,\end{equation}
where the minimal value is:
\begin{equation}\label{argmin_c2uneVit}
c_{min}^2=4\,{\frac { \left( -\rho^+\rho^-{(c^-)}^{2}+(\rho^+)^{2}
{(c^+)}^{2}+(\rho^-)^{2}{(c^-)}^{2}-\rho^-\rho^+{c_{{1
}}}^{2} \right) {(c^+)}^{2}\rho^-{(c^-)}^{2}\rho^+}{(\rho^
{{+}})^{4}{(c^+)}^{4}-2\,(\rho^+)^{2}{(c^+)}^{2}(\rho^-)^{
2}{(c^-)}^{2}+(\rho^-)^{4}{(c^-)}^{4}}}\,.
\end{equation}
Thus, for water and air at $25^{\circ}\text{C} $ it is found that $c_{min}\approx 23.78 \,m/s\,$. Plots of $c_w$ is shown in Figures \ref{fig1j} to \ref{fig1j1}.
Note, if it is supposed that $\frac{\rho^+}{\rho^-} \ll  \frac{c^-}{c^+} = O(1)$, the following approximations are obtained:
\begin{equation}\label{argmin_c2uneVit_dl}
c_{min}^2 \approx 4\,{\left( \frac{\rho^+}{\rho^-}\right) (c^+)^2}\,,
\end{equation}
\begin{equation}\label{argmin_alphauneVit_dl}
\alpha_{min}^2 \approx \,{\frac{1}{2}\left( 1 + \frac{\rho^+}{\rho^-} \frac{(c^-)^2 - (c^+)^2}{(c^+)^2}\right) }\,.
\end{equation}
So the minimum speed of sound is approximately $ 2 c^+\sqrt{\frac{\rho^+}{\rho^-}}$, which is attained around the value $\alpha^+ = {1}/{2}$. 

It is also noted that Wood's speed of sound corresponds 
to the following thermodynamic derivative of the mixture equation of state \eqref{isoS1u4_bis}:
\begin{equation}\label{c_wood_thermo1}
c^2 = \left. \frac{\partial p}{\partial \rho} \right|_{\eta}\,.
\end{equation}

To obtain the mixture EOS \eqref{isoS1u4_bis}, we note that given each fluid equation of state $p=P(\rho^{\pm})$, the mixture with an additional state variable $\alpha^+$ or $\eta$ and the mechanical equilibrium constraint $p^+ =p^-$\,, is a divariant system with an equation of state that can be expressed in the form $p = P(\rho,\alpha^+)$ or $p = P(\rho,\eta)$.

\subsubsection{A model with two temperatures}

As noted before, for smooth solutions this system is equivalent to 
\begin{equation}\label{1b80_bis}
\rho_t+div(\rho\, u)=0\,,
\end{equation}
\begin{equation}\label{2b80_bis}
(\rho u)_t+div(\rho u\otimes u)+ \nabla p=0\,,
\end{equation}
\begin{equation}\label{3b80_bis}
\eta_t+ u.\nabla\eta =0\,,
\end{equation}
\begin{equation}\label{4b80_bis}
s^{\pm}_t+ u.\nabla s^{\pm} =0\,,
\end{equation}
\begin{equation}\label{6b80_bis}
    p = P(\rho,\eta,s^+,s^-)\,.
\end{equation}

The eigenvalues for the linearized system are $u+c$, $u-c$, $u$ (with multiplicity $3$), with $c$ given by Wood's formula:
\begin{equation}\label{c_x_bis}
\frac{1}{\rho^2 c^2}=\frac{\eta^+}{(\rho^+)^2 (c^+)^2} +\frac{\eta^-}{(\rho^-)^2 (c^-)^2} \,,
\end{equation}
and with the definition $P = P(\rho,s^+,s^-,\eta)$: 
\begin{equation}\label{c_x_ter}
c^2 = \left. \frac{\partial P}{\partial \rho} \right|_{s^+\,,s^-\,,\eta}\,.
\end{equation}

\paragraph{Link to thermodynamics \ref{AppenA} }
Here, given each fluid equation of state $P=P(\rho_i, s_i)$, the mixture (in mechanical equilibrium $p^+ =p^-$, but not necessary in thermal equilibrium) is a quadri-variant system with an equation of state of the form $P = P(\rho, s^+, s^-, \eta)$. This is referred to in Temkin \cite{Temkin92} as the frozen equilibrium equation of state. Also, Stewart and Wendroff \cite{Stewart} note that a two-fluid flow is an out of equilibrium process and that Wood's formula is an out of equilibrium speed of sound. The reader is referred to Landau and Lifshitz \cite[p. 42]{Landau2} for a discussion on the characteristic times for mechanical and thermal equilibriums between two phases. A complete equilibrium speed of sound can be defined with the additional hypothesis $T^+ = T^-$, as shown in the next Section. 



\subsubsection{A model with thermal equilibrium}

Linearizing the system (\ref{1u1T1}-\ref{1u1T4}) or its equivalent form (\ref{1u1T1_bis}-\ref{1u1T5_bis}), leads to a set of eigenvalues of the form $\{u-c, u, u,u+c\}$. The equation of state in (\ref{1u1T5_bis}) is the same as the one defined in \ref{thermo:equilibrium} and after diagonalization of the system, it was found that the speed of sound $c$ is given by:
$$
c^2 =\left.\frac{\partial p}{\partial \rho}\right|_{\eta\,,s}\,,
$$
which is the same as the thermodynamic speed of sound for a two-fluid mixture (\ref{cS-equilibrium}) and is referred to in Temkin \cite{Temkin92} as the relaxed equilibrium speed of sound.

\subsection{Two velocity models}

\subsubsection{Isentropic case}\label{sec521}

By setting $W = ~^t( p \,, \alpha^+ \,, u^+ \,, u^-)$, it can be shown that system (\ref{2u1}) - (\ref{2u2}) is equivalent for smooth solutions to :
\begin{multline}
\label{eqn_p}
\left[\alpha^+\rho^-(c^-)^{2}+\alpha^-\rho^+(c^+)^{2}\right]\frac{\partial p}{\partial t} +(\alpha^+\rho^-(c^-)^{2}u^++\alpha^-\rho^+(c^+)^{2}u^-)\nabla p \\ + \rho^+\rho^-(c^+)^{2}(c^-)^{2} \text{div}(\alpha^+u^+ +\alpha^-u^-)=0\,,
\end{multline}
\begin{multline}
\label{eqn_alpha}
\left[\alpha^+\rho^-(c^-)^{2}+\alpha^-\rho^+(c^+)^{2}\right]\frac{\partial
\alpha^+}{\partial t}+\alpha^+\alpha^-(u^+-u^-)\nabla p \\
+\alpha^-\rho^+(c^+)^{2} \text{div}(\alpha^+u^+)- \alpha^+\rho^-(c^-)^{2} \text{div}(\alpha^-u^-)=0\,,
\end{multline}
\begin{equation}
\label{eqn_u_2}
\frac{\partial u^\pm}{\partial t}
+\frac{1}{\rho^\pm}\frac{\partial p}{\partial x} +u^\pm\frac{\partial u^\pm}{\partial x} =0\,,
\end{equation}
\begin{equation}\label{2u3_tier}
EoS_\pm(p,\rho^\pm)=0\,.
\end{equation}

Note that (\ref{eqn_p}) to (\ref{eqn_u_2}) can be written as a quasilinear system~:
\begin{equation}
\label{quasi_w}
\frac{\partial W_i}{\partial t}+ \sum^3_{k = 1} \sum^{4}_{j = 1}B_{i,j}^k(W)\frac{\partial W_j}{\partial x_k} = 0\,.
\end{equation}
Let ${W_0}$ be a given constant state. The linearization of (\ref{quasi_w}) around this state is~:
\begin{equation}
\label{lquasi_w}
\frac{\partial W_i}{\partial t}+ \sum^3_{k = 1} \sum^{4}_{j = 1}B_{i,j}^k({W_0})\frac{\partial W_j}{\partial x_k} = 0\,.
\end{equation}

For the case of a constant solution $W_0$ with $u^{-} = u^{+} = u_0$, and after diagonalization of the matrices $B^k$, it is found that all the eigenvalues are in the set $ \{u_0, u_0 \pm c_0\}$ with $c_0$ the dynamic speed of sound:
\begin{equation}\label{vitson2vit}
c_{0} = \sqrt{{\frac {{(c^-)}^{2}{(c^+)}^{2} \left( \alpha^+\rho^-+
\alpha^-\rho^+ \right) }{\alpha^+\rho^-{(c^-)}^{2}+
\alpha^-\rho^+{(c^+)}^{2}}}}\,.
\end{equation}
The index $0$ in $c_0$ refers to the fact there is ``a zero force'' in the momentum equations that acts to equalize the two $u^+$ and $u^-$ velocities\footnote{Like {\it e.g.} the added mass force (\ref{f1}).}.\\
Observe that $c_0$ is a monotone function (homographic) of the volume fractions $\alpha^+$ or $\alpha^-$. This can also be seen in Figures \ref{fig1j} and \ref{fig1j1}.

Considering only the equations of state (\ref{2u3_tier}), the `thermodynamic speed of sound', which can be seen as a compressibility coefficient: 
$$
c^2 =\left. \frac{\partial p}{\partial \rho} \right|_{\eta}\,,
$$
is still given by Wood's formula (\ref{c_x_bis}), and hence, in this case, it is different from the actual propagation speed of sound (\ref{vitson2vit}) in the medium.

\subsubsection{A model with energy}
The energy equations (\ref{2u4}) can be rewritten for the entropies $s^\pm$ as:
\begin{equation}\label{seqn}
 s^\pm_t + u^\pm \nabla s^\pm = 0\,,
\end{equation}
and equations (\ref{eqn_p})-(\ref{eqn_u_2}) remain the same as in the isentropic case. In view of (\ref{seqn}) and given the block structure of the Jacobian matrix, only the multiplicity of the eigenvalue $u_0$ is altered, while the speed of sound is the same as in the isentropic case:

\begin{equation}\label{vitson2vit_e}
c_0^2 ={\frac  {{(c^-)}^{2}{(c^+)}^{2} \left( \alpha^+\rho^-+
\alpha^-\rho^+ \right) }{\alpha^+\rho^-{(c^-)}^{2}+
\alpha^-\rho^+{(c^+)}^{2}}}\,.
\end{equation}

\subsubsection{A model with added mass term}\label{sec323}

Following the calculation in the Appendix of Redford {\it et al.} \cite{Redford}, linearization of System (\ref{1.8}) to (\ref{1.12}) leads to the following speed of sound:

\begin{equation}\label{SoS1}
c_\kappa^2\equiv \left(\alpha^+\frac{\rho+\kappa\rho^+}{\rho(1+\kappa)\rho^+}+
\alpha^-\frac{\rho+\kappa\rho^-}{\rho(1+\kappa)\rho^-}\right)\frac{\rho^+\rho^-(c^+)^2(c^-)^2}{\alpha^+\rho^-(c^-)^2+\alpha
^-\rho^+(c^+)^2}\,.
\end{equation}
Then one can easily compute:
\begin{equation}\label{SoS2}
c_0^2\equiv \frac{(\alpha^-\rho^++\alpha^+\rho^-)(c^+)^2(c^-)^2}{\alpha^+\rho^-(c^-)^2+\alpha
^-\rho^+(c^+)^2}\,,
\end{equation}
and
\begin{equation}\label{SoS3}
c_\infty^2\equiv \frac{\rho^+\rho^-(c^+)^2(c^-)^2}{\rho(\alpha^+\rho^-(c^-)^2+\alpha
^-\rho^+(c^+)^2)}\,.
\end{equation}
where $\rho=\alpha^+\rho^++\alpha^-\rho^-$. Note that (\ref{SoS2}) is the same as  (\ref{vitson2vit}), whereas (\ref{SoS3}) is Wood's formula (\ref{SoS1}) for the one-velocity model: $c_\infty=c_w$\,.\\

It can easily be seen that:
\begin{equation}\label{ordre-c_0_c_w}
c_0 \ge c_w\quad  \mbox{ and }\quad  c_0 = c_w \quad  \mbox{ if and only if }\quad \alpha^+\,\alpha^-\,(\rho^+-\rho^-)^2=0\,.
\end{equation}

As expected, $c_\kappa$ interpolates between $c_0$\,, the speed of sound for the two-fluid two-velocity model in Section \ref{sec521}, and $c_w\,,$  the speed of sound for the two-fluid one-velocity model in Section \ref{sec311}. Moreover, $c_\kappa^2$ is a simple linear interpolation between $c_0^2$ and $c_w^2$:
\begin{equation}\label{interpol}
c_\kappa^2=\frac{\kappa\,c_w^2+c_0^2}{\kappa +1}\,.
\end{equation}
Furthermore, thanks to $c_0 \ge c_w$\,, $c_\kappa$ is a strictly convex interpolation between $c_0$ and $c_w$ (as long as $\alpha^+\,\alpha^-\,(\rho^+-\rho^-)^2>0$)\,.


\subsubsection{Baer-Nunziato Model for Reactive Granular Materials}\label{BNModel}
The 1986 article of Baer and Nunziato \cite{Baer_Nunziato} introduced a two-phase model describing flame spread and deflagration-to-detonation in gas-permeable reactive granular materials.
\paragraph{Convective part of Baer-Nunziato Model}
This reads as:
\begin{equation}\label{BN1}
    \alpha^+_t + u^+ \alpha^+_x=0\,,
\end{equation}
\begin{equation}\label{BN2}
    (\alpha^+\rho^+)_t + (\alpha^+\rho^+u^+)_x=0\,,
\end{equation}
\begin{equation}\label{BN3}
    (\alpha^+\rho^+u^+)_t + (\alpha^+(\rho^+(u^+)^2+p^+))_x
    -p^-\alpha^+_x=0\,,
\end{equation}
\begin{equation}\label{BN4}
    (\alpha^+\rho^+E^+)_t + (\alpha^+(\rho^+E^++p^+)u^+)_x
    -p^-u^+\alpha^+_x=0\,,
\end{equation}
\begin{equation}\label{BN5}
    (\alpha^-\rho^-)_t + (\alpha^-\rho^+u^-)_x=0\,,
\end{equation}
\begin{equation}\label{BN6}
    (\alpha^-\rho^-u^-)_t + (\alpha^-(\rho^-(u^-)^2+p^-))_x
    +p^-\alpha^+_x=0\,,
\end{equation}
\begin{equation}\label{BN7}
    (\alpha^-\rho^-E^-)_t + (\alpha^-(\rho^-E^-+p^-)u^-)_x
    +p^-u^+\alpha^+_x=0\,,
\end{equation}
where subscript $+$ (resp. $-$) refers to the solid (resp. gas). This system is closed by two equations of state (EoS):
\begin{equation}\label{BN8}
    EoS_\pm(p^\pm\,,\rho^\pm\,,e^\pm)=0\,.
\end{equation}
Although this model is not a two-fluid model, but rather a two-phase (gas-solid) model, it is included in this paper because some authors (this subject will be touched upon in Section \ref{Sec5}) use it in simulation codes for both two-phase and two-fluid flows. As will now be seen, the values for the speed of sound are particular in this case.
\paragraph{The Bar-Nunziato model Speed of Sound}\label{Sec3.4}
According to Embid and Baer \cite{Embid}, the convection operator eigenvalues that correspond to the Bar-Nunziato model (\ref{BN1})-(\ref{BN8}) are as follows:
\begin{equation}
    u^+-c^+\,,u^++c^+\,,u^+ \mbox{(twice)}, u^--c^-\,,u^-+c^-\,,u^-\,.
\end{equation}
Using the notation in (\ref{plane_wave1}), here:
\begin{equation}\label{LambdaBM}
\Lambda=\{-c^+\,,-c^-\,,0\,,c^+\,,c^-\}\,\mbox{ (recall that }u^\pm=0\,). \end{equation}
Hence, the speeds of sound for this model are the phasic speeds of sound:
\begin{equation}
    c^+ \quad\mbox{ and }\quad c^-\,.
\end{equation}
\section{Dispersion relations}\label{Sec4}
The first goal in this Section is to derive the dispersion relation for the three kinds of models: $1$ fluid, $2$ fluids with $1$ velocity and $2$ fluids with $2$ velocities. Then various asymptotic relations are derived for the speed and attenuation of sound-waves (sometimes called pressure-waves in the literature).
The results are summarized in Table \ref{tab5} (Section \ref{Sec4.4}), where the expressions for these quantities are the same provided that relevant values are taken for the speed of sound at infinite length scale ($k=0$)\footnote{Or vanishing frequency $\omega=0\,.$} and for the equivalent dynamic viscosity. 
Concerning the single-fluid case, reference is made to the classical works of Stokes \cite{Stokes}, Kirchhoff \cite{Kirchhoff} and Fletcher \cite{Fletcher} and asymptotic developments are presented with the specific physical conditions under which the different dispersion and attenuation formulas apply.
\subsection{A single conductive and viscous fluid}\label{AppenNS}
The goal here is to illustrate both the nature of a dispersion relation in fluid mechanics and the technique used to derive it. Considering the 1D compressible Navier-Stokes equations: 
\begin{equation}
\left\{ 
\begin{array}{l}
 \dfrac{\partial \rho}{\partial t}+ u \dfrac{\partial \rho}{\partial x} + \rho \dfrac{\partial u}{\partial x}=0\,,\vspace{0.2cm}\\
 \dfrac{\partial u }{\partial t}+ u \dfrac{\partial u}{\partial x}
+\dfrac{1}{\rho} \dfrac{\partial p}{\partial x} =\dfrac{1}{\rho}\dfrac{\partial}{\partial x}\left(\dfrac{4\,\mu}{3}\dfrac{\partial u}{\partial x}\right)\,,\vspace{0.2cm}\\
\dfrac{\partial s}{\partial t}+ u \dfrac{\partial s}{\partial x} = \dfrac{4\,\mu}{3\,\rho\, T} \left(\dfrac{\partial u}{\partial x}\right)^2 + \dfrac{1}{\rho\, T}  \dfrac{\partial}{\partial x}\left(\lambda\,\dfrac{\partial T}{\partial x}\right)\,,\vspace{0.2cm}\label{systeme_NST} 
\end{array}
\right.
\end{equation}
where dynamic viscosity $\mu$\,, thermal conductivity $\lambda$ and pressure $p$ are known as a function of density and temperature. Using Stokes' hypothesis for the second viscosity coefficient $\mu'=\zeta +2\,\mu/3=0$, where $\zeta$ is the second Lamé's elastic coefficient for the fluid (usually denoted by $\lambda$ which is kept here for the thermal conductivity). In all the formulas hereafter, the general case $\mu' \geq 0$ can be retrieved by replacing $\mu$ by $\mu_{eq} = \mu + 3\,\mu'/4 = 3\,\mu/2 + 3\,\zeta/4\, $, and redefining the Prandt number with $\mu_{eq}$ accordingly.

The linearization of \eqref{systeme_NST} around a constant solution $\rho_0, u_0=0, T_0\,$ reads:
 
\begin{equation}
\left\{ 
\begin{array}{l}
 \dfrac{\partial \rho}{\partial t}+ \rho_0 \dfrac{\partial u}{\partial x}=0\,,\vspace{0.2cm}\\
 \dfrac{\partial u }{\partial t}
+\dfrac{1}{\rho_0} \dfrac{\partial p}{\partial x} =\dfrac{4\,\mu_0}{3\,\rho_0}\dfrac{\partial^2 u}{\partial x^2}\,,\vspace{0.2cm}\\
\dfrac{\partial s}{\partial t}=\dfrac{\lambda_0}{\rho_0\,T_0} \dfrac{\partial^2 T}{\partial x^2} = \dfrac{\lambda_0\Gamma_0}{\rho_0^2}\dfrac{\partial^2 \rho}{\partial x^2}+ \dfrac{\lambda_0}{\rho_0\,C_{v,0}}\dfrac{\partial^2 s}{\partial x^2}\,,\vspace{0.2cm}\\
p= c_0^2 \rho + \rho_0\,\Gamma_0 T_0 \,s\,,
\label{systeme_NST_L} 
\end{array}
\right.
\end{equation}
where (see Table \ref{tab4} for notation):
\begin{equation}\label{Cv_etc} 
\Gamma\equiv \frac1\rho \frac{\partial p}{\partial e} \big|_{\rho}\,,\quad c^2 \equiv \frac{\partial p}{\partial \rho} \big|_s\,,\quad
C_{v}\equiv \frac{\partial e}{\partial T} \big|_{\rho}=T\,\frac{\partial s}{\partial T} \big|_{\rho}\,.
\end{equation}
The second equality in the $s$ evolution equation of \eqref{systeme_NST_L} follows from the two thermodynamic identities \eqref{Relation_8}.
In general the linearization should be made around a constant solution $\rho_0\,,$ $ u_0\,,$ $T_0\,,$ but using the Galilean invariance of the Navier-Stokes equations it can be assumed that $u_0=0$\,. 
Hence, when $u_0\neq0$ all the eigenvalues should be shifted by $u_0$.\\

The linear differential system (\ref{systeme_NST_L}) is of the form $\frac{\partial W}{\partial t} + A \frac{\partial W}{\partial x}= B \frac{\partial^2 W}{\partial x ^2}$, where $W = ~ ^t(\rho\,,u\,, s)$ and: 
$$
A = \left( \begin{array}{ccc}
0 & \rho_0 & 0 \\
\frac{c_0^2}{\rho_0} & 0 & \Gamma_0\,T_0  \\ 
0 & 0 & 0 \end{array} \right),\quad
B = 
\left( \begin{array}{ccc}
0 & 0 & 0 \\
0 & \frac{4\,\mu_0}{3\,\rho_0} & 0  \\ 
\frac{\lambda_0\Gamma_0}{\rho_0^2} & 0 & \frac{\lambda_0}{\rho_0\,C_{v,0}} \end{array} \right).
$$
Non-vanishing plane-wave solutions of the form \eqref{plane_wave_bis} are of interest. Hence, ${\omega_R}/{k}$ represents the propagation speed of the plane-wave and the dispersion relation corresponds to the 
roots 
of the characteristic polynomial of  the $3\times 3$ matrix: 
  $-i\omega I + i k A + k^2 B\,. $
Dropping the 0 subscript, this gives:
\begin{eqnarray}\label{RD_NS}
    \left(\frac{\omega}{k}\right)^{3} +   \frac{i\,k }{\rho} \left(\frac{4\,\mu}{3} + \frac{\lambda}{C_v}\right)\left(\frac{\omega}{k}\right)^{2}-  \left(c^{2} + \frac{4\,k^{2}\, \lambda\, \mu}{3\,\rho^{2}\,C_v }\right)\left(\frac{\omega}{k}\right)-\frac{i\, k\,  \lambda\,c_T^2}{ \rho\,C_v}\,,
\end{eqnarray}
using the identity $c_T^2= c^2 - \Gamma^2 C_v T $\,, see \eqref{Relation_4}.\\

Using the Prandtl number and Stokes' attenuation factor $a$\,, \eqref{CT}, equation \eqref{RD_NS} can be rewritten as \eqref{NS1D}, which will now be analysed.
\begin{itemize}
    \item[(i)] For $\mu= 0$ and $\lambda=0$ (or $a_1=0$, $\text{\text{Pr}}=\infty$), the speed of sound ${\omega}/{k} \in \{ 0, \pm c  \}$ is recovered for the inviscid Euler equation. This system is hyperbolic and non dispersive. 
    \item[(ii)] In the case of a 
    non heat conducting ($\lambda=0$\,, or $\text{Pr}=\infty$) viscous flow  ($\mu \neq 0$), Stokes' attenuation and dispersion relations are recovered \cite{Stokes}: 
 \begin{equation}\label{RDStokes}
 \left(\frac{\omega}{k}\right)^{3}  +  \frac{4\,i\, k\,\mu}{3\,\rho} \left(\frac{\omega}{k}\right)^{2} -  c^{2}\,\left(\frac{\omega}{k}\right) = 0\,,\end{equation}
  \begin{equation}\label{RDStokes1}
  \Rightarrow  \frac{\omega}{k} \in \left\{ 0, {- i \frac{2\,k\, \mu}{3\rho} } \pm  { \sqrt{ c^2 - \frac{4\,k^2 \mu^2}{9\,\rho^2} } }\right\}.    \end{equation}
Let us consider a plane wave \eqref{plane_wave}. Here $1/k$ represents a characteristic length for the wave. In the context of continuum mechanics (as opposed to rarefied flows) where the Navier-Stokes equations are valid, the Knudsen number $\text{Kn}$ built with this characteristic length should not be greater than the critical Knudsen number $\text{Kn}_c=10^{-2}$\,, that is:
\begin{equation}\label{Knudsen}
    \text{Kn}\equiv\frac{k\,\mu}{\rho\,c}\le \text{Kn}_c=10^{-2} \,.
\end{equation}
Hence:
\begin{eqnarray}\nonumber
  c(k) 
&=&c\sqrt{ 1 - \frac{4\,(\text{Kn})^2}{9} } = c\,\left(1- \frac{2\,(\text{Kn})^2}{9} + \mathcal{O}(\text{Kn}^4)\right)\\\label{dev_Knud}
&\approx&c\,\left(1- \frac{2\,k^2 \mu^2}{9\,\rho^2\,c^2}\right)=c-\frac{2\,k^2 \mu^2}{9\,\rho^2\,c}\,.
\end{eqnarray}

 According to \eqref{plane_wave_bis}, the disturbances are a linear combination of functions of the form:
 \begin{equation}\label{plane_wave_NS}
    W^{\pm}= W_0\,\exp{-\frac{2\,k^2\,\mu\,t}{3\,\rho}}\exp{ik\,(x\pm c(k)\,t)}\,,
\end{equation}
\begin{equation}\label{plane_wave_NS0}
    W^{0}= W_0\,\exp{-\frac{2\,k^2\,\mu\,t}{3\,\rho}}\exp{ik\,x}\,.
\end{equation}
The waves \eqref{plane_wave_NS} are dispersive since their velocities are $\pm c(k)$\,, where:
\begin{equation}\label{c_de_k}c(k)\equiv { \sqrt{ c^2 - \frac{4\,k^2 \mu^2}{9\,\rho^2} } }= \sqrt{ c^2 - k^2\,a_1^2}
\approx c-\frac{2\,k^2 \mu^2}{9\,\rho^2\,c} \,,
\end{equation}
depends on $k$\,. Moreover they are exponentially damped with time.

    \item[(iii)] For $\mu=0$ and $\lambda\neq 0$ ($\text{Pr}=0$) the dispersion relation \eqref{RD_NS} reads:
    \begin{equation}
    \left(\frac{\omega}{k}\right)^{3} +   \frac{i\,k\,\lambda }{\rho\,C_v} \left(\frac{\omega}{k}\right)^{2}-  c^{2} \left(\frac{\omega}{k}\right)-
 \frac{i\, k\,  \lambda\,c^2}{\gamma\, \rho\,C_v} = 0\,.
   \label{GV1}
    \end{equation}
In Vekstein \cite{Vekstein} this dispersion relation is given in the case of a perfect gas. Introducing the Knudsen number based on thermal diffusivity:
\begin{equation}\label{GV2}
    \text{Kn}_{th}\equiv \frac{k\,\lambda}{\rho\,C_p\,c}\,,
\end{equation}
then \eqref{GV1} can be rewritten as:
 \begin{equation}\label{GV3}
    \left(\frac{\omega}{k}\right)^{3} +   i\,\gamma\,c\,\text{Kn}_{th} \left(\frac{\omega}{k}\right)^{2}-  c^{2} \left(\frac{\omega}{k}\right)-i\,c^3\,\text{Kn}_{th} = 0\,.
    \end{equation}
In the context of fluid dynamics the Knudsen number is small and the techniques used in \ref{AppenC} lead immediately to the following asymptotic expressions for the three roots, $(\omega/k)^\pm = (\omega_R^\pm + i \,\omega_I^\pm)/k  $ and $(\omega/k)^0 = (\omega_R^0 + i\, \omega_I^0)/k$, of (\ref{NS1D}).
\begin{equation}\label{GV4}
    \frac{\omega^\pm_I}{k}= -\frac{(\gamma-1)\,k\,\lambda}{2\,\rho\,\gamma\,C_v}+\mathcal{O}\left(\frac{k^2\,\lambda^2}{\rho^2\,C_v^2\,c}\right)\,,\quad \frac{\omega^\pm_R}{k}= \pm c+\mathcal{O}\left(\frac{k^2\,\lambda^2}{\rho^2\,C_v^2\,c}\right)\,,
\end{equation}        
        \begin{equation}\label{GV5}
    \frac{\omega_I^0}{k}= -\frac{k\,\lambda}{\rho\,\gamma\,C_v}+\mathcal{O}\left(\frac{k^2\,\lambda^2}{\rho^2\,C_v^2\,c}\right)\,,\quad \frac{\omega_R^0}{k}=\mathcal{O}\left(\frac{k^2\,\lambda^2}{\rho^2\,C_v^2\,c}\right)\,,
\end{equation}
and \eqref{GV4} generalizes Vekstein \cite{Vekstein} to the case of arbitrary divariant fluids. Note also that the zero order term in $\frac{\omega^\pm_I}{k}$ is the same as the thermal contribution in the Stokes-Kirchhoff \cite{Kirchhoff} attenuation: $$\frac{(\gamma-1)\,k\,\lambda}{2\,\rho\,C_p} = \frac{(\gamma-1) \text{Kn}_{th}}{2} c = \frac{\lambda \, k}{2\,\rho\,C_v} \left( 1- \frac{c_T^2}{c^2}\right).$$
\end{itemize}

For the general case of a viscous and heat conducting fluid, the reader is referred to \cite{BG} for an analysis of the dispersion relation \eqref{RD_NS}.
\subsection{A one-velocity two-fluid viscous model}\label{sec6.2}
The system reads:
\begin{equation}\label{1u1visc01}
\left\{ 
\begin{array}{l}
(\alpha^\pm\rho^\pm)_t+div(\alpha^\pm\rho^\pm u)=0\,,\vspace{0.2cm}\\
(\rho u)_t+div(\rho u\otimes u)+ \nabla p= \text{div}({\bm\sigma}), \vspace{0.2cm}\\
(\alpha^\pm\rho^\pm E^\pm)_t+div(\alpha^\pm\rho^\pm H^\pm u)+
 p\,\alpha^\pm_t=  \text{div}(\alpha^\pm \,{\bm\sigma}^\pm\cdot u)\,,
\end{array}
\right.
\end{equation}
where $p = P^\pm(\rho^\pm,s^\pm)$ and with Stokes hypothesis : 
\begin{equation}\label{1u1visc1}
 {\bm\sigma}^\pm=\mu^\pm\,(\nabla u+^t\nabla u)-\frac{2\,\mu^\pm}{3}\,div\,u\,{\bf I}\,,\quad {\bm\sigma}=\mu\,(\nabla u+^t\nabla u)-\frac{2\,\mu}{3}\,div\,u\,{\bf I} \,.
\end{equation}

It is more convenient to rewrite the system in the non-conservative form with the variables  $\eta$, $\rho$, $u$ and $s^{\pm}$. In one dimension we obtain :
\begin{equation}
\left\{ 
\begin{array}{l}
 \dfrac{\partial \eta}{\partial t}+ u \dfrac{\partial \eta}{\partial x} =0,\vspace{0.2cm}\\
 \dfrac{\partial \rho}{\partial t}+ u \dfrac{\partial \rho}{\partial x} + \rho \dfrac{\partial u}{\partial x}=0,\vspace{0.2cm}\\
 \dfrac{\partial u }{\partial t}+ u \dfrac{\partial u}{\partial x}
+\dfrac{1}{\rho} \dfrac{\partial p}{\partial x} =\dfrac{1}{\rho}\dfrac{\partial }{\partial x}\left(\dfrac{4\,\mu}{3} \dfrac{\partial u}{\partial x}\right),\vspace{0.2cm}\\
\dfrac{\partial s^{\pm}}{\partial t}+ u \dfrac{\partial s^{\pm}}{\partial x} = \dfrac{4\,\alpha^{\pm}\mu^{\pm}}{3\,\rho^\pm T^\pm} \left(\dfrac{\partial u}{\partial x}\right)^2,\vspace{0.2cm}\\
p = P^\pm(\rho^\pm,s^\pm)\,.\label{systeme+vitNS} 
\end{array}
\right.
\end{equation}

As for the monofluid case, the second viscosity can be taken into account by replacing $\mu_\pm$ by $\mu_{eq}^{\pm} = 3\,\mu_\pm/2 + 3\,\zeta_\pm/4\, $, where $\zeta_\pm$ are the second Lamé coefficients for the two-fluids and redefining the average viscosity $\mu_{eq}$ accordingly by $\alpha^+\,\mu_{eq}^++\alpha^-\,\mu_{eq}^-$.

The obtained dispersion equation is 
$$ \left(\frac{\omega}{k}\right)^{5} +  \frac{4\,i\mu\,k }{3\,\rho}  \left(\frac{\omega}{k}\right)^{4}  - c_w^{2}  \left(\frac{\omega}{k}\right)^{3}  = 0\,. $$
where $c_w$ is Wood's speed of sound,
$
c_w^2 = \left.\frac{\partial p}{\partial \rho}\right|_{\eta\,,s^\pm} 
$,
given by formula (\ref{c_x}).

The dispersion formula is the same as for the monofluid case, using the average quantities $\mu$, $\rho$ and Wood's speed of sound $c_w$\,: 
$$
c(k) = \frac{\omega}{k} \in \left\{ 0, {- i \frac{2\,k\, \mu}{3\rho} } \pm  { \sqrt{ c_w^2 - \frac{4\,k^2 \mu^2}{9\,\rho^2} } }\right\}.   $$
Hence, introducing the {\it ad hoc} Knudsen number
\begin{equation}\label{Knudsen_w}
    \text{Kn}_w\equiv\frac{k\,\mu}{\rho\,c_w}\le \text{Kn}_c=10^{-2} \,,
\end{equation}
and then, using the same approach as for the single fluid case, for the sound-waves (with notation in \eqref{plane_wave_bis}-\eqref{def_c(k)}):
\begin{equation}\label{c_de_k_Wood}c(k)= { \sqrt{ c_w^2 - \frac{4\,k^2 \mu^2}{9\,\rho^2} } }
= c_w-\frac{2\,k^2 \mu^2}{9\,\rho^2\,c_w}+\mathcal{O}\left(\frac{2\,k^4 \mu^4}{\rho^4\,c_w^3}\right) \,,\quad \sigma(k)=\frac{2\,k^2\,\mu}{3\,\rho}\,.
\end{equation}





\begin{figure}[ht!]
\begin{center}
\includegraphics[scale=0.6]{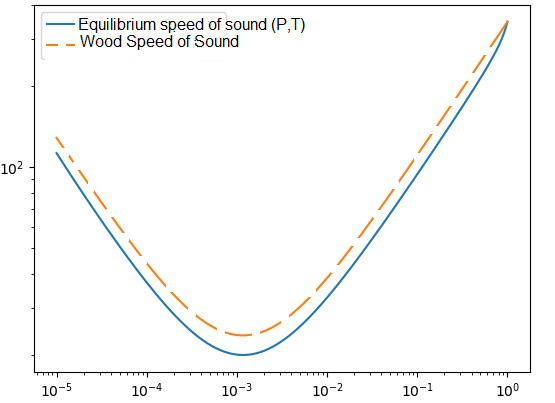}
\caption{Plot of $c_w$ and $c_{eq}$ for a water-air mixture (foam) against the water mass fraction on the horizontal axis. The thermal equilibrium leads to a slightly smaller speed of sound.} \label{fig_w_wq}
\end{center}
\end{figure}

\noindent {\it Thermal conductivity: } If one considers thermal conductivity  $\lambda^\pm$ for each fluid, the entropy equations become: 
$$
\frac{\partial s^\pm}{\partial t}+ u \frac{\partial s^\pm}{\partial x} = \frac{4\,\alpha^\pm \,\mu^\pm}{3\,\rho^\pm\, T^\pm} \left(\frac{\partial u}{\partial x}\right)^2 + \frac{1}{\rho^\pm\, T^\pm}  \frac{\partial}{\partial x}\left(\lambda^\pm\,\frac{\partial T^\pm}{\partial x}\right)\,,
$$
and after linearization:
$$
\frac{\partial s^\pm}{\partial t}=\frac{\lambda^\pm}{\rho^\pm\,T^\pm} \frac{\partial^2 T^\pm}{\partial x^2} = \frac{\lambda^\pm\Gamma^\pm}{\rho^\pm {}^2}\frac{\partial^2 \rho^\pm}{\partial x^2}+ \frac{\lambda^\pm}{\rho^\pm\,C_{v,\pm}}\frac{\partial^2 s^\pm}{\partial x^2}\,,\vspace{0.2cm}
$$
$$
p= c_w^2 \rho + \rho^+\,\Gamma^+ T^+ \,s^+ + \rho^-\,\Gamma^- T^- \,s^-\,.
$$
Calculations similar to the monofluid case (Section \ref{AppenNS}) can be conducted to derive Stokes-Kirchoff type dispersion and attenuation relations and their asymptotic analysis. Here, only thermal transfers within each phase are taken into account. A more elaborate model would also take into account the thermal exchanges between the two phases, through their interfaces.\\
Another thermal model consists of considering the one-temperature model (\ref{1u1T1})-(\ref{1u1T4}) with an effective conductivity for the two-fluid mixture. The dispersion and attenuation from the monofluid case applies exactly, where $c$ is replaced by $c_{eq}$ in \eqref{cS-equilibrium}. It is recalled that $c_{eq}$ is different to Wood's velocity, although they can be numerically close, as seen in Figure \ref{fig_w_wq}. 

\subsection{A two-velocity two-fluid viscous model, with added mass}\label{sec6.3}
Following Ishii and Hibiki \cite{Ishii}, viscosity is taken into account in the system (\ref{1.8})-(\ref{1.12}). The mass conservation equations are, of course, unchanged while the right-hand side parts of (\ref{1.10}) and (\ref{1.12}) are respectively modified as follows:
\begin{equation}
\label{1.10v}
\mbox{r.h.s. viscous momentum eqn. }
= p\nabla \alpha^\pm+f^\pm+ \text{div}(\alpha^\pm{\bm \sigma^\pm})\,,
\end{equation}
\begin{equation}
\label{1.12v}
\mbox{r.h.s. viscous energy eqn. } = f^\pm\cdot u^\pm+ \text{div}(\alpha^\pm\mu^\pm\,{\bm \sigma^\pm}\cdot u^\pm)\,.
\end{equation}
The dispersion relation for this model is derived in \ref{AppenE} and it is found that:
\begin{multline}\label{RD2vit}
\left(\frac{\omega}{k}\right)^{3}  +
i k \left(  \frac{ 4\mu_{\kappa} }{3\rho \left(\kappa + 1\right) } \right)
\left(\frac{\omega}{k}\right)^{2}\,
- \left ( c_{\kappa}^2  + \frac{ 16\mu^+ \mu^- k^{2}  }{9(\kappa+1) \rho^+ \rho^-} \right)
\left(\frac{\omega}{k}\right) \\
- \frac{4i k \overline{\mu}\,\rho\,c_w^2}{3\rho^+\,\rho^-\,(\kappa + 1)}  = 0\,,
\end{multline}

where:
\begin{equation}\label{mukappa}
\mu_\kappa\equiv (1+\kappa)\,\mu+\alpha^+\mu^-\frac{\rho^+}{\rho^-}+\alpha^-\mu^+\frac{\rho^-}{\rho^+} = \kappa\,\mu + \rho \left(\frac{\mu^+}{\rho^+} + \frac{\mu^-}{\rho^-} \right)\,,
\end{equation}

\begin{equation}\label{barmugamma}
\quad\overline{\mu}\equiv \alpha^+\mu^-+\alpha^-\mu^+ \,.
\end{equation}
\begin{Remarque}
In the case of vanishing viscosities $\mu^\pm=0$, as expected, the dispersion relation (\ref{RD2vit}) leads to (\ref{SoS1}).
\end{Remarque}

\paragraph{Asymptotic speed of sound for second order $\mu$}
Introducing the phasic Knudsen numbers:
\begin{equation}\label{Knudsen_pm}
\text{Kn}^\pm\equiv\frac{k\,\mu^\pm}{\rho^\pm\,c^\pm}\le \text{Kn}_c=10^{-2} \,,
\end{equation}
where the $\omega/k$ solutions of the dispersion relation \eqref{RD2vit} that corresponds to sound-waves can be expanded with first order $\text{Kn}^\pm$ as follows:
\begin{equation}\label{dev2vit1}
\frac{\omega^\pm}{k} = \pm c_\kappa + i\,{\bf\tilde{c}}_+\,\text{Kn}^+ + i\,{\bf\tilde{c}}_-\,\text{Kn}^- +\mathcal{O}((\text{Kn}^+)^2)+(\text{Kn}^-)^2))\,,
\end{equation}
where:
\begin{equation}\label{dev2vit2}
{\bf\tilde{c}}_\pm\equiv {\frac {2 \left( {\rho}^{2}\alpha^\mp c_w^2-\rho^\mp  ( \rho +\kappa\alpha^\pm\rho^\pm) c_\kappa^2 \right)c{{^\pm }} }{ 3\left(
\kappa+1 \right)\,\rho\, \rho^\mp c_\kappa^2}}\,.
\end{equation}
If this expansion is expressed in terms of $\mu^\pm$ then:
\begin{equation}\label{dev2vit3}
\frac{\omega^\pm}{k} = \pm c_\kappa -\frac{2\, i\,k}{3\,\rho}\,({\bm \varphi}_+\,\mu^+ +{\bm \varphi}_-\,\mu^-) +\mathcal{O}\left(\left(\frac{k\,\mu^+}{\rho^+\,c^+}\right)^2+\left(\frac{k\,\mu^-}{\rho^-\,c^-}\right)^2\right)\,,
\end{equation}
\begin{equation}\label{dev2vit4}
{\bm \varphi}_\pm\equiv {\frac {  \rho^\mp  ( \rho +\kappa\alpha^\pm\rho^\pm)c_\kappa^2  -{\rho}^{2}\alpha^\mp c_w^2   }{ \left(
\kappa+1 \right)\,\rho^+\,\rho^- c_\kappa^2}}\,.
\end{equation}
Moreover, the second order $\mu^\pm$ term in \eqref{dev2vit3} is as follows:
\begin{equation}\label{dev2vit5}
-\frac{2\,k^2}{9\,\rho^2\,c_\kappa}({\bm \theta}^\pm_{+}(\mu^+)^2+2\,{\bm \delta}^\pm\,\mu^+\,\mu^-+{\bm \theta}^\pm_{-}(\mu^-)^2)
 +\mathcal{O}\left(\left(\frac{k\,\mu^+}{\rho^+\,c^+}\right)^3+\left(\frac{k\,\mu^-}{\rho^-\,c^-}\right)^3\right)\,,
\end{equation}
\begin{equation}\label{dev2vit6}
{\bm \theta}^\pm_\ell\equiv \pm{\frac {{\bm \varphi}_\ell \left( (\kappa+1)\rho^-
\rho^+{\bm \varphi}_\ell\,c_\kappa^{2}+4\,\rho^{2}\alpha^{-\ell} c_w^{2}
 \right) }{(\kappa+1)\rho^+\rho^-c_\kappa^{2}}}\,,
 \quad \ell\in\{+,-\}
\,,
\end{equation}
\begin{equation}\label{dev2vit7}
{\bm \delta}^\pm\equiv\pm {\frac {c_\kappa^{2}(\kappa+1)\rho^-\rho^+{\bm \varphi}_{{+}}
{\bm \varphi}_{{-}}+2\,c_{{w}}^{2}{\rho}^{2}(\alpha^-\,{\bm \varphi}a_{{-}}+\alpha^+{\bm \varphi}_{{+}})-2\,c_{{\kappa}
}^{2}{\rho}^{2}}{(\kappa+1)\rho^+\rho^-c_\kappa^{2}}}
\,.
\end{equation}








\subsection{Summary for viscous models}\label{Sec4.4}
The results given in Table \ref{tab5} are obtained if pressure-waves for the $3$ previous viscous models are considered.
\begin{table}[ht!]
    \centering
    \begin{tabular}{|c|c|c|}\hline 
       Model  & Speed of sound $c(k)$ & Attenuation $\sigma(k)$ \\\hline\hline
       & & \\
      $1$ fluid, $1$ velocity  & $c-\frac{2\,k^2 \mu^2}{9\,\rho^2\,c}$ & $\frac{2\,k^2\,\mu}{3\,\rho}$ \\
      Section \ref{AppenNS}& & \\\hline
      & & \\
      $2$ fluids, $1$ velocity & $c_w-\frac{2\,k^2 \mu^2}{9\,\rho^2\,c_w}$ & $\frac{2\,k^2\,\mu}{3\,\rho}$ \\
      Section \ref{sec6.2}& & \\\hline
      & & \\
      $2$ fluids, $2$ velocities & $c_\kappa-\frac{2\,k^2 \mu_\theta^2}{9\,\rho^2\,c_\kappa}$ & $\frac{2\,k^2\,\mu_\varphi}{3\,\rho}$ \\
     Section \ref{sec6.3} & & \\\hline
    \end{tabular}
    \caption{Asymptotic speed of sound and attenuation for second order Knudsen numbers and viscous models, $\mu_\varphi\equiv{\bm \varphi}_+\,\mu^+ +{\bm \varphi}_-\,\mu^-$\,, $\mu_\theta^2\equiv{\bm \theta}^+_{+}(\mu^+)^2+2\,{\bm \delta}^+\,\mu^+\,\mu^-+{\bm \theta}^+_{-}(\mu^-)^2$\,, see \eqref{dev2vit4}, \eqref{dev2vit6} and \eqref{dev2vit7}.}
    \label{tab5}
\end{table}

\section{Consequences for CFD methods}\label{Sec5}
As already discussed, there is no universal model for two-fluid or two-phase flows, in contrast with single-fluid flow where the Navier-Stokes equation is the reference model for dense (as opposed to rarefied) flows. A CFD model is the discrete version of a continuous mathematical model. Hence, the mathematical models will now be discussed.\\

Roughly speaking, as far as averaged models are concerned, there are two principle classes of models: one and two velocity models. In each of these two classes, a chosen equilibrium assumption can be made based on the physical situation ({\it e.g.} single temperature assumption). Two velocity models are more general in the sense that one can often recover a one-velocity model by adding the two conservation of momentum equations together. 
For example, summing the two equations (\ref{2u2}) leads to (\ref{isoS1u2}).\\

Hence, let us discuss two velocity models. In addition to the classical transport term of the form $(\alpha^\pm\rho^\pm u^\pm)_t+div(\alpha^\pm\rho^\pm u^\pm\otimes
u^\pm)$, these models have two momentum equations, which involve pressure terms, viscous terms and terms that model the exchange of momentum between the two fluids. Since these terms are obtained from the phasic Navier-Stokes equations by an averaging process, Ishii and Hibiki \cite{Ishii}, they must be modeled in order to obtain a closed system. For example, these may include the drag force between the two fluids or the added mass term. These modeling terms will of course depend strongly on the particular flow under investigation, including the regimes concerned.\\

Most modern CFD codes rely on Finite Volume Methods (FVM) in order to locally satisfy the balance equations. For these methods, the propagation direction of the waves in the system is critical. Indeed, FVM applied to CFD are based on the construction of flux terms at an edge between two volumes. These fluxes are determined by examining the propagation direction of the information (or waves).
As discussed in the Introduction, it is therefore of the utmost importance that the considered model capture physically relevant speeds of sound. When dealing with the flow of two compressible fluids, as supported by experiments (see Section \ref{sec142}), one expects that the physical model on which the code is built will propagate pressure-waves in a physically meaningful manner. 
That is, at a single speed that depends on the constant physical state that is perturbed. In other words, physical models that produce two sound speeds ({\it e.g.} one per fluid) are not suitable.\\

In contrast, the model of Baer and Nunziato \cite{Baer_Nunziato} has two sound speeds: one for the gas phase and the other for the solid phase. But here, the Authors are dealing with gas-permeable, reactive granular materials and not two-fluid compressible flows. From the physical point of view this could be interpreted as follows for large solid volume fraction. In such a case the solid grains may touch and a pressure wave can propagate from one solid particle to another, hence producing sound waves at the solid sound speed. See {\it e.g.} Mézière {\it et al.} \cite{Meziere}.\\

The first two-fluid CFD codes were developed in the context of thermal-hydraulics in the 1970s (see Lyczkowski \cite{Lyczkowski}). This initial effort was mainly focused on physical models because the subject was new at this time. Concerning numerical methods, the approach was essentially the so-called ``pressure based methods", which consists of deriving an elliptic equation for the pressure. The advantage of this kind of method is that it naturally handles low Mach number flows ({\it i.e} quasi-incompressible flows). However, these methods are overly diffusive such that the solution gradients are smeared, in particular the pressure gradients which are essential for sizing of installations. Moreover, the method is implicit, which allows to solve large and ill-conditioned linear systems.\\

As far as Finite Volumes are concerned, it is only since the 1990s that convection methods ({\it i.e} use of the pressure waves in the computation of the fluxes) have been used. The great advantage of convection methods is that they allow accurate computation of strong solution gradients. Most convection methods are explicit, leading to ``simple" codes. However, use of these methods is sometimes challenging for low Mach number flows. This drawback can be corrected using preconditioning techniques, see {\it e.g.} Guillard and Viozat \cite{Guillard}.\\

Most modern industrial codes are based on convection methods. Hence correctly computing the pressure-waves is an important issue. Some Authors (Saurel and Abgrall \cite{Saurel}, Gallouet {\it et al.} \cite{Gallouet}, Tian {\it et al.} \cite{Tian}, Städtke {\it et al.} \cite{Stadtke}) propose two-fluid flow models with inconsistent sound speeds. Namely, the two distinct speeds of sound of each fluid, as for the Baer and Nunziato model. In fact, these authors are led to this because they need to compute exact solutions (use of Riemann invariant or Riemann solvers) in order to apply a robust numerical method, which is derived for a single fluid flow, such as Godunov's method \cite{Godunov} or Roe's scheme \cite{Roe}, and is then applied to multi-fluid flow. This is not always possible for the physical models
, but in these articles, the model is chosen so that the numerical method can be applied. Moreover, in the aforementioned references, mostly $1D$ validations are proposed and the current authors are unaware of applications of such methods to industrial cases.\\

In Ghidaglia {\it et al.} \cite{GKL_CRAS, GKL}, a characteristic flux finite volume method has been proposed. It was designed especially for multi-fluid or multi-phase models. This numerical method is applicable here regardless of the convective part of the model. The method is based on the reduction of the convection matrix ({\it e.g.} $A$ in \eqref{A_1} or $B^k$ in \eqref{lquasi_w}) for the considered model. As shown in Ghidaglia \cite{FLUX}, this method is a natural extension of Roe's scheme, that is derived for single fluid flows, to multi-fluid flows. This method has been used by various authors for $1D$, $2D$ and $3D$ simulations and for industrial applications in the context of multi-fluid and multi-phase flows (Halama {\it et al.} \cite{Halama}, Kervella {\it et al.} \cite{Kervella}, Dutykh {\it et al.}  \cite{DUTYKH}, Sahmim {\it et al.} \cite{Sahmim}, Boucker \cite{Boucker}, Rovarch \cite{Rovarch}, Benjelloun {\it et al.} \cite{Benjelloun}, Redford {\it et al.} \cite{Redford}, ...).

\section{Conclusions and perspectives}\label{Sec6}
Various expressions have been derived for the speed of sound for two-fluid models by taking convective and viscous effects into account. The consequences of the model choice for CFD codes have been discussed.\\
In further work, the important case of two fluids of the same substance in thermodynamic equilibrium with phase change will be addressed. The effect of surface tension on the value of the speed of sound is also to be considered in the case of two fluids.\\

Another question that needs careful investigation in the future is the following. Considering a numerical code designed to simulate a given mathematical model, for which there is an analytical expression for the speed of sound and the attenuation. 
Is there consistency and convergence of values from the numerical model with respect to those found with the analytical expression?

\section*{Acknowledgement}
The authors would like to thank John Redford and Valentin Leroy for the careful reading of this paper and for their valuable remarks for the improvement of this paper. 

\appendix

\section{The thermodynamic speed of sound}\label{AppenA}

In thermodynamics, for single-phase fluids and mixtures, it is a common abuse of language to refer to the `speed of sound' for any quantity of the form 
\begin{equation}\label{eq:thermo_c2}
c^2 = \frac{\partial p}{\partial \rho}\,,
\end{equation}
where $\rho$ and $p$ are respectively the density and the pressure of the system defined as averages or in some other way.

A monophase fluid is a divariant thermodynamical system, and the equation (\ref{eq:thermo_c2}) implies that a thermodynamic state variable is held constant. Entropy and temperature are the two quantities usually considered. The isentropic speed of sound, that is still denoted as $c$, is the velocity observed in normal physical circumstances which arises from linear perturbation theory for the Euler system :
\begin{equation}\label{eq:thermo_c2_s}
c^2 = \left. \frac{\partial p}{\partial \rho}\right|_s\,.
\end{equation}

The non conductive speed of sound, that is denoted here as $c_T$, is reported in some experiments as the actual speed of propagation of sound under some physical circumstances (boundary layers, microfluidic, high frequency...)\cite{Blackstock,Fletcher}, and this is justified for physical reasons:
\begin{equation}\label{eq:thermo_c2_T}
c_T^2 = \left. \frac{\partial p}{\partial \rho} \right|_T\,.
\end{equation}

These definitions for the non conductive and isentropic speeds of sound can be generalized for two phase mixtures, that are or are not in thermodynamic equilibrium. 

\subsection{Equilibrium Speed of Sound in two phase mixtures}\label{thermo:equilibrium}

A two-fluid mixture (with no mass exchange between the two fluids), with mass fraction $\eta^+$ of fluid $+$ and $\eta^-$ for fluid $-$, is under complete thermodynamic equilibrium if $T^+ = T^-$ and $P^+ = P^-$. 
Such a system, which is assumed to be constantly in equilibrium, is also a divariant system and the mixture thermodynamic quantities, $v\,,\rho\,,e\,,s$ given in Table \ref{tab3},
satisfy the classical thermodynamic identities ({\it e.g.} $de = T\,ds - p\,dv$... etc). Hence, the different  thermodynamic coefficients can be defined from partial derivatives, such as the isentropic and non conductive speed of sound. Furthermore, all the thermodynamic relations (Reech, Mayer, ... etc) are verified by these coefficients. In particular we recall the following identities that will be used hereafter:
\begin{equation}
c_T^2 = c^2 - \Gamma^2 T C_v\,,
\end{equation}
and
\begin{equation}
\label{CfromCT}
\frac{1}{\rho^2 c_s^2} = \frac{1}{\rho^2 c_T^2} - T \frac{\chi^2}{C_p}\,,
\end{equation}
where $\Gamma$ is the Grüneisen parameter defined as $\Gamma=v\frac{\partial p}{\partial e}]_v$, $C_p$ is the isobaric heat capacity, and $\chi$ is the coefficient of the isobar thermal expansion, which is defined as $\chi = \frac{1}{v}\frac{\partial v}{\partial T}]_p$\,. 

For such system, it can be proven that the equilibrium non-conductive speed of sound is given by Wood's formula:
\begin{equation}\label{cT-equilibrium}
\frac{1}{\rho^2 c_T^2} = \frac{\eta^+}{(\rho^+)^2 (c_{T}^+)^2} + \frac{\eta^-}{(\rho^-)^2 (c_{T}^-)^2}\,,
\end{equation}
or equivalently
\begin{equation}
\label{cT-equilibrium-bis}
\frac{1}{\rho c_T^2} = \frac{\alpha^+}{\rho^+ (c_{T}^+)^2} + \frac{\alpha^-}{\rho^- (c_{T}^-)^2}\,. \end{equation}
Indeed, one can derive the mixture volume $v$ in Table \ref{tab3} with respect to $P=P^+ =P^-$ at constant $T=T^+=T^-$, to obtain (\ref{cT-equilibrium}).

The isentropic speed of sound can then be computed using the relation (\ref{CfromCT}) to give

\begin{equation}
\label{cS-equilibrium}
\frac{1}{\rho^2 c^2} = \frac{\eta^+}{(\rho^+)^2 (c_{T}^+)^2} + \frac{\eta^-}{(\rho^-)^2 (c_{T}^-)^2} - T \frac{\chi^2}{\eta C_p^+ + (1-\eta) C_p^-}\,,
\end{equation}
with $\chi$ being the mixture isobar thermal expansion given by $$\chi = \frac{1}{v}\frac{\partial v}{\partial T}]_p = \frac{1}{v} (\eta^+\chi^+ v^++ \eta^-\,\chi^- v^-) =  \alpha^+\, \chi^+ +  \alpha^-\, \chi^-\,,$$
and the mixture isobaric heat capacity $T \frac{\partial s}{\partial T} |_P$ is given by: 
\begin{equation}
C_p = \eta^+\, C_p^+ + \eta^-\,C_p^-\,.
\end{equation}

The speed of sound (\ref{cS-equilibrium}) is referenced in \cite{Temkin92,Temkin00} as the relaxed equilibrium speed of sound and a different, although equivalent, formula is given in \cite{Temkin92,Temkin00}:
$$
\frac{ (c^-)^2}{c^2} =  \frac{\rho}{\rho^-} \left( \gamma^-\,\alpha^- + \gamma^+ \alpha^+ \frac{ \rho^-(c^-)^2}{\rho^+(c^+)^2}\right) - (\gamma^- -1) \frac{\alpha^- + \alpha^+\, \chi^+ /\chi^-}{ \eta^- + \eta^+\, C_p^+/C_p^-}.
$$

Note that, for intermediate values of $\alpha^\pm$, this formula gives lower values than Wood's formula, hence it corrects Wood's overestimation as seen in Figure \ref{figkarplus}. Moreover, as the third term in $(\ref{cS-equilibrium})$ is dominated by the gas terms $( (\chi^-)^2 >> (\chi^+)^2)$ this correction gives similar results to equation (\ref{c_woodT}), which is shown in Figure \ref{figkarplus}. \\   

\subsection{Non-equilibrium speed of sound in two phase mixtures}

The equilibrium speed of sound presented in the previous two sections is meant to be applied to systems in complete thermodynamic equilibrium, but in many situations this hypothesis is not justified. Particularly in the context of fluid flows, where the assumption of thermal equilibrium between the two phases does not hold \cite{Stewart}. However, mechanical equilibrium (i.e. the equality of pressure) is a more rapid phenomenon and can be assumed in most cases \cite[page~44]{Landau2}. In this context and assuming a mixture without exchange of matter between the two phases, the speed of sound is given by Wood's equation \cite{kieffer,Brennen}:

$$
c^2 =\left. \frac{\partial \rho}{\partial P}\right|_{s^+,s^-,\eta}\,,
$$
$$
\frac{1}{\rho c^2} = \frac{\alpha^+}{\rho^+ (c^+)^2}  + \frac{\alpha^-}{\rho^- (c^-)^2}\,.  
$$

Given the above justification, this speed is referred to in the literature \cite{Temkin92,Stewart} as the frozen equilibrium speed of sound. 

\section{On the dispersion relation for a single fluid 
}\label{AppenC}
In the two first Sections of this Appendix our goal is to study the dispersion relation \eqref{RD_NS} that has been rewritten as $P\left(\frac{\omega}{k}\right)=0$\, (see \eqref{NS1D}), where (referring to \eqref{CT} and Table \ref{tab4} for notation):
\begin{equation}\label{A3.1}
P(X)\equiv X^{3} + 2\,i\,k\,a_1  \,\left( 1+\frac{3\,\gamma}{4\,Pr} \right)X^{2}  -  \left(c^{2} + \frac{3\,\gamma}{\text{Pr}}a_1^2\,k^2\right)X - i\,\frac{3\,k\,a_1\,c^2}{2\,Pr}\,,
\end{equation}
where the last term has been modified using the identity $c_T^2= {c^2}/{\gamma}$\,, see \eqref{Relation_4}.\\
Introducing the two second-degree polynomials $Q$ and $Q_T$\,:
\begin{eqnarray}\label{A3.2}
Q(X)\equiv X^{2}  + 2\,i\,k\,a_1  \, X -  c^{2}\,,\\
\label{A3.3}
Q_T(X)\equiv X^{2}  + 2\,i\,k\,a_1  \, X -  {c^2}/{\gamma}\,,
\end{eqnarray}
it is seen that (note that these two polynomials are independent of $Pr$):
\begin{equation}\label{A3.4}
P(X) = X\,Q(X)+i\,k\,a_1\,\frac{3\,\gamma}{2\,Pr}\,Q_T(X)\,.
\end{equation}
\subsection{Asymptotic for large Prandtl number}\label{AppenC.1}
When the Prandtl number $Pr$ is infinite, which is the case for a 
non-conductive ($\lambda=0$) viscous fluid  ($\mu \neq 0$), $P(X)=X\,Q(X)$ and \eqref{RDStokes} is recovered.\\
Since the three $P$ roots are distinct in this case, one can easily see\footnote{Using the implicit function theorem, as is done hereafter in the proof of Proposition~\ref{prop20}.} that for large Prandtl number ({\it i.e.} for $\lambda \ll C_p\,\mu$), that is the case where viscous diffusion is more prevalent 
than thermal diffusion, the three $P$ roots have a smooth dependence on $1/Pr$ and can be expanded as:
\begin{equation}\label{A3.5}
X_{\text{Pr}}^{(\ell)}=X_0^{(\ell)}-i\,k\,a_1\,\frac{3\,\gamma}{2\,Pr}
\frac{Q_T(X_0^{(\ell)})}{Q(X_0^{(\ell)})+X_0^{(\ell)}\,Q'(X_0^{(\ell)})}+\mathcal{O}\left(\frac{1}{Pr^2}\right)\,,
\end{equation}
where $X_0^{(\ell)}$ are the three $X\,Q(X)$ roots. According to \eqref{RDStokes1}, with $c(k)$ defined in \eqref{c_de_k}, these roots are:
\begin{equation}\label{A3.6}
X_0^{(-1)}=- i k\,a_1 -  c(k)\,,\quad X_0^{(0)}=0\,,\quad X_0^{(+1)}=-ik\,a_1 +  c(k)\,,
\end{equation}
and after some computations, the following result is deduced from \eqref{A3.5}.
\begin{prop}\label{prop15}
For large Prandtl number, the three $P$ roots satisfy:
\begin{eqnarray}\label{A3.7}
X_{\text{Pr}}^{(\mp)}=-ik\,a_1 \mp  c(k)
\mp\frac{3\,(\gamma-1)\,k\,a_1\,(k\,a_1\pm ic(k))}{4\,c(k)}\frac{1}{\text{Pr}}+\mathcal{O}\left(\frac{1}{Pr^2}\right)\,,\\
\label{A3.8}
X_{\text{Pr}}^{(0)}=-i\,\frac{3\,k\,a_1}{2\,Pr}+\mathcal{O}\left(\frac{1}{Pr^2}\right)\,.
\end{eqnarray}
\end{prop}
\subsection{Asymptotic for small Prandtl number}\label{AppenC.2}
Let us now address the case where the Prandtl number is small ({\it i.e.} for $C_p\,\mu \ll \lambda$), which is the case where thermal diffusion is more prominent with respect to viscous diffusion.\\
At the limit where $\text{Pr}=0$, the dispersion relation $P\left(\frac{\omega}{k}\right)=0$ reads $Q_T\left(\frac{\omega}{k}\right)=0$ but $Q_T$ has only two roots while $P$ has three. It will be proven that for $Pr\ll1$ the two roots of $P$ situated are on curves starting from the two roots of $Q_T$ (as in \eqref{A3.5}) and the third root is large $\mathcal{O}\left(\frac{1}{\text{Pr}}\right)$, see Propositions \ref{prop10} and \ref{prop20}.\\
Indeed, denoting the three roots of $P$ by $\xi^{(\ell)}_{\text{Pr}}$\,, where $\ell=-1\,,0$ and $1$, and if $\xi_{\text{Pr}}^{(- 1)}$ (resp. $\xi_{\text{Pr}}^{(+1)}$) is close to $\xi^{(- 1)}$ (resp. $\xi^{(+1)}$), where $\xi^{(\ell)}$ are the roots of $Q_T$, that is:
\begin{equation}\label{A3.11}
\xi^{(\ell)}\equiv - i k\,a_1 +\ell \, c_T(k)\,,\quad c_T(k)\equiv \sqrt{ c_T^2 - k^2\,a_1^2}\,,\quad \ell=\pm 1\,,
\end{equation}
which will now be proven in the following.
\begin{prop}\label{prop10}
The following asymptotic expressions are obtained:
\begin{eqnarray}\label{A9.1}
\lim_{Pr\rightarrow 0}\xi_{\text{Pr}}^{(\mp)}=- i k\,a_1 \mp \, c_T(k)\,,\\
\label{A9.2}
\xi_{\text{Pr}}^{(0)}\sim i\frac{3\,k\,a_1\,\gamma}{2\,Pr}\,,\mbox{ as } {Pr\rightarrow 0}\,
.
\end{eqnarray}
\end{prop}
This result is readily derived from a simple observation.
\begin{lem}
The three $\xi^{(\ell)}_{\text{Pr}}$ roots of $P$ satisfy the identity:
\begin{equation}\label{A3.10}
\xi_{\text{Pr}}^{(-1)}\,\xi_{\text{Pr}}^{(0)}\,\xi_{\text{Pr}}^{(+1)}=-i\frac{3\,k\,a_1\,c^2}{2\,Pr}\,.
\end{equation}
\end{lem}
This relation is obvious because according to \eqref{A3.4}, $P(0)=i\,k\,a_1\,\frac{3\,\gamma}{2\,Pr}\,Q_T(0)$, and then Proposition \ref{prop10} follows.\\
Here again, Proposition \ref{prop10} can be refined in the spirit of Proposition \ref{prop15}, and the following result is proven.
\begin{prop}\label{prop20}
For small Prandtl numbers, the three $P$ roots satisfy:
\begin{eqnarray}\label{A3.12}
\xi_{\text{Pr}}^{(\mp)}=- i k\,a_1 \mp \, c_T(k)\pm \frac{(\gamma-1)\,c^2\,(k\,a_1\mp i\,c_T(k))}{3\,\gamma^2\,k\,a_1\,c_T(k)}Pr+\mathcal{O}(Pr^2)\,,\\
\label{A3.13}
\xi_{\text{Pr}}^{(0)}= i\frac{3\,k\,a_1\,\gamma}{2\,Pr}+ \mathcal{O}(Pr)\,
.
\end{eqnarray}
\end{prop}
\paragraph{Proof} $\xi_{\text{Pr}}^{(\ell)}$ are the solutions of $P(X)=0$ or equivalently:
\begin{equation}\label{A3.15}
F(X\,,Pr)\equiv Pr\,X\,Q(X)+i\,\frac{3\,k\,a_1\,\gamma}{2}\,Q_T(X)=0\,.
\end{equation}
The function $F$ is smooth and $\frac{\partial F}{\partial X}(X\,,0)=i\,\frac{3\,k\,a_1\,\gamma}{2}\,Q'_T(X)$\,. For each $\ell=\pm 1$, $F(\xi^{(\ell)}\,,0)=i\,\frac{3\,k\,a_1\,\gamma}{2}\,Q_T(\xi^{(\ell)})=0$ and since the roots $\xi^{(\ell)}$ of $Q_T$ are simple: $\frac{\partial F}{\partial X}(\xi^{(\ell)}\,,0)\neq 0$. Hence, two smooth curves are found to exist using the implicit function theorem for small $Pr$, such that $F(\xi_{\text{Pr}}^{(\ell)}\,,Pr)=0$ and $\xi_{0}^{(\ell)}=\xi^{(\ell)}$\,. Then \eqref{A3.12} follows immediately from the first-order Taylor expansion of $F$ with respect to $X$ and $Pr$ at the point $(\xi^{(\ell)}\,,0)\,.$

Concerning  \eqref{A3.13}, \eqref{A3.12} is used together with the identity \eqref{A3.10}.
\section{The dispersion relation for a two-velocity two-fluid viscous model with added mass}\label{AppenE}
The system \eqref{1.8}-\eqref{f1}, where the right hand sides of (\ref{1.10}) and (\ref{1.12}) are respectively modified according to (\ref{1.10v}) and (\ref{1.12v}), is addressed here.
One dimensional flows can only be considered as far as the dispersion relation is concerned. Indeed, plane-waves $ W= W_0\,\exp{i(k\,\cdot x-\omega\,t)}$\,, $x\in\mathbb{R}^3$\,, $k\in\mathbb{R}^3$\,, propagate in the direction of the  wave vector $k$, and when $k\neq 0$ the system satisfied by $\omega$ and $k$ can be found by writing the corresponding $1D$ system, where $x\in\mathbb{R}$ is the direction of $k$\,.
By choosing variables $W = ~^t(  \alpha^+ \,, p \,, u^+ \,, u^-\,,s^+\,,s^-)$, this system is equivalent (for smooth solutions) to :

\begin{multline}
(\alpha^- \rho^+ (c^+)^2 + \alpha^+ \rho^- (c^-)^2 )\frac{\partial \alpha^+}{\partial t} + \alpha^+ \alpha^- (u^+ - u^- ) \frac{\partial p}{\partial x}  
+\alpha^+ \alpha^- (c^+)^2 \rho^+ \frac{\partial u^+}{\partial x} \\ - \alpha^+ \alpha^- (c^-)^2 \rho^- \frac{\partial u^-}{\partial x} 
+ \left[ \alpha^- \rho^+ u^+ (c^+)^2 + \alpha^+ \rho^- u^- (c^-)^2 \right] \frac{\partial \alpha^+}{\partial x} \\ =  \alpha^+ \alpha^-  \frac{4\,\mu^+}{3}\left(\frac{\partial u^+}{\partial x}\right)^2  - \alpha^+\alpha^-  \frac{4\,\mu^-}{3}\left(\frac{\partial u^-}{\partial x}\right)^2\,,
\end{multline}

\begin{multline}
( \alpha^+ \rho^-  (c^-)^2  + \alpha^- \rho^+ (c^+)^2)  \frac{\partial p}{\partial t} + (\alpha^+ \rho^- (c^-)^2  u^+ + \alpha^- \rho^+ (c^+)^2 u^- ) \frac{\partial p}{\partial x}+  \\
+\alpha^+ \rho^+ \rho^- (c^-)^2  (c^+)^2  \frac{\partial u^+}{\partial x} + \alpha^-  \rho^+ \rho^- (c^+)^2 (c^-)^2 \frac{\partial u^-}{\partial x}  
 + \left[ \rho^+ \rho^-  (c^-)^2 (c^+)^2\right] (u^+ - u^-)  \frac{\partial \alpha^+}{\partial x} \\ =  \alpha^+ \rho^- (c^-)^2   \frac{4\,\mu^+\,\Gamma^+}{3}\left(\frac{\partial u^+}{\partial x}\right)^2 +  \alpha^- \rho^+ (c^+)^2    \frac{4\,\mu^-\,\Gamma^-}{3}\left(\frac{\partial u^-}{\partial x}\right)^2\,,
\end{multline}

\begin{multline}
\frac{\partial u^+}{\partial t} +\frac{\rho+\kappa\rho^+}{\rho(1+\kappa)\rho^+}\frac{\partial p}{\partial x}
+\frac{\rho+\kappa\alpha^+\rho^+}{\rho(1+\kappa)}u^+\frac{\partial u^+}{\partial x}+\frac{\kappa\alpha^-\rho^-}{\rho(1+\kappa)}u^-\frac{\partial u^-}{\partial x} = 
\\
\frac{\rho+\kappa\alpha^+\rho^+}{\rho(1+\kappa)}  \frac{\mu^+}{ \alpha^+\rho^+}  \frac{\partial \alpha^+}{\partial x} \frac{\partial u^+}{\partial x} 
+
 \frac{\rho+\kappa\alpha^+\rho^+}{\rho(1+\kappa)}\frac{\mu^+}{\rho^+} \frac{\partial^2 u^+}{\partial x^2} +
\\
+\frac{\kappa \mu^- }{\rho(1+\kappa)} \frac{\partial \alpha^-}{\partial x} \frac{\partial u^-}{\partial x} 
+\frac{\kappa\alpha^- \mu^-}{\rho(1+\kappa)} \frac{\partial^2 u^-}{\partial x^2}\,,
\end{multline}

\begin{multline}
\frac{\partial u^-}{\partial t}+\frac{\rho+\kappa\rho^-}{\rho(1+\kappa)\rho^-}\frac{\partial p}{\partial x}
+\frac{\kappa\alpha^+\rho^+}{\rho(1+\kappa)}u^+\frac{\partial u^+}{\partial x}+\frac{\rho+\kappa\alpha^-\rho^-}{\rho(1+\kappa)}u^-\frac{\partial u^-}{\partial x}
= 
\\
\frac{\kappa\mu^+}{\rho(1+\kappa)}  \frac{\partial \alpha^+}{\partial x} \frac{\partial u^+}{\partial x} + \frac{\kappa\alpha^+ \mu^+}{\rho(1+\kappa)} \frac{\partial^2 u^+}{\partial x^2}  +
\\
+\frac{\rho+\kappa\alpha^-\rho^-}{\rho(1+\kappa)} \frac{\mu^-}{ \alpha^-\rho^-}  \frac{\partial \alpha^-}{\partial x} \frac{\partial u^-}{\partial x} +
\frac{\rho+\kappa\alpha^-\rho^-}{\rho(1+\kappa)}\frac{\mu^-}{\rho^-} \frac{\partial^2 u^-}{\partial x^2}  \,,
\end{multline}

\begin{equation}
s^+_t+u^+ \frac{\partial s^+}{\partial x}= \frac{4\,\mu^+}{3\,\rho^+ T^+} \left(\frac{\partial u^+}{\partial x}\right)^2 \,,
\end{equation}
\begin{equation}
s^-_t+u^- \frac{\partial s^-}{\partial x}= \frac{4\,\mu^-}{3\,\rho^- T^-} \left(\frac{\partial u^-}{\partial x}\right)^2 \,,
\end{equation}

\begin{equation}
p=P(\rho, \eta, s^+, s^-) \,.
\end{equation}

Keeping only the linear terms, it is seen that the entropy equations can be decoupled. The dispersion relation is then given when the matrix determinant vanishes: 
$$ -i\omega I + i k A + k^2 B\,,$$ 
with :
\begin{eqnarray} 
& A(\alpha^+, p, u^+, u^-)= & \nonumber \\
&\left [
\begin {array}{cccc} 
\frac{\alpha^- \rho^+ (c^+)^2 u^++\alpha^+ \rho^- (c^-)^2 u^-}{\overline{\pi}} & \frac{\alpha^+ \alpha^-\left (u^+-u^-\right)}{\overline{\pi}} &\frac{(\alpha^+ \alpha^- \rho^+ (c^+)^2)}{\overline{\pi}}& -\frac{(\alpha^+\alpha^-\rho^- (c^-)^2)}{\overline{\pi}}\\
\frac{\rho^+\rho^-(c^+)^2(c^-)^2 \left (u^+-u^-\right) }{\overline{\pi}}&\frac{\alpha^+\rho^-(c^-)^2u^++\alpha^-\rho^+ (c^+)^2u^-}{\overline{\pi}} &\frac{\alpha^+\rho^+\rho^-(c^+)^2(c^-)^2}{\overline{\pi}}&\frac{\alpha^-\rho^+\rho^-(c^+)^2(c^-)^2}{\overline{\pi}}\\
0&\frac{\rho+\kappa\rho^+}{\rho(1+\kappa)\rho^+}& \frac{\rho+\kappa\alpha^+\rho^+}{\rho(1+\kappa)}u^+&\frac{\kappa\alpha^-\rho^-}{\rho(1+\kappa)}u^-  \\
0&\frac{\rho+\kappa\rho^-}{\rho(1+\kappa)\rho^-}& \frac{\kappa\alpha^+\rho^+}{\rho(1+\kappa)}u^+ & \frac{\rho+\kappa\alpha^-\rho^-}{\rho(1+\kappa)}u^- \end{array}
\right ] &\,\nonumber
\end{eqnarray}
$$
\mbox{ where }\quad\overline{\pi} \equiv  \alpha^+ \rho^-  (c^-)^2  + \alpha^- \rho^+ (c^+)^2=\frac{\rho^+\rho^-(c^+)^2(c^-)^2}{\rho\,c_w^2}\,,
$$

$$
B = \left[
 \begin {array}{cccc}
 0& 0 & 0& 0\\
 0& 0 & 0& 0 \\
 0& 0 & \frac{\rho+\kappa\alpha^+\rho^+}{\rho(1+\kappa)}\frac{4\mu^+}{3\rho^+}& \frac{4\kappa\alpha^- \mu^-}{3\rho(1+\kappa)} \\
  0& 0 & \frac{4\kappa\alpha^+ \mu^+}{3\rho(1+\kappa)} & \frac{\rho+\kappa\alpha^-\rho^-}{\rho(1+\kappa)}\frac{4\mu^-}{3\rho^-}
  \end{array}
\right].
$$

Then the determinant, for $u^+ =u^- =0$, vanishes according to:
\begin{multline}
\left(\frac{\omega}{k}\right)^{4}  +
4i k \left(  \kappa \frac{ \alpha^+  \mu^+ + \alpha^- \mu^- }{3\rho \left(\kappa + 1\right)} + \frac{\mu^+ \rho^- + \mu^- \rho^+}{3\rho^+ \rho^- \left(\kappa + 1\right)} \right)
\left(\frac{\omega}{k}\right)^{3} \\
- \left ( k^{2} \frac{16 \mu^+ \mu^-  }{9(\kappa+1) \rho^+ \rho^-}
+ (c^+)^2 (c^-)^2 \frac{ \kappa \rho^+ \rho^- 
+  \rho (\alpha^+ \rho^- + \alpha^-  \rho^+)
} {{\overline{\pi}} (\kappa+1)\rho }\right)
\left(\frac{\omega}{k}\right)^2 \\
- \frac{4i k (c^+)^2 (c^-)^2 \left( \alpha^- \mu^+ 
+ \alpha^+ \mu^-  
\right)}{3{\overline{\pi}} (\kappa + 1)}\left(\frac{\omega}{k}\right)=0\,.
\end{multline}

By noting that $c_{\kappa}^2 =\frac{ \kappa \rho^+ \rho^- 
+  \rho (\alpha^+ \rho^- + \alpha^-  \rho^+)
} {{\overline{\pi}} (\kappa+1)\rho^+ \rho^- }\,c_w^2$, the dispersion relation can be rewritten, eliminating the first vanishing root $\left(\frac{\omega}{k}\right) =0$\,, as with \eqref{RD2vit}.

\section{Some classical thermodynamic identities used in this article}\label{AppenD}

The dispersion relations derived in this paper involve some thermodynamic coefficients that depend on the equations of state. All the coefficients are well referenced in the literature, but the notation is not universal. Furthermore, they are expressed in the case of a perfect gas in many references. For example, in this case the Grüneisen coefficient $\Gamma$ is equal to $\gamma-1$, where $\gamma$ is the ratio of the heat capacity at constant pressure to the heat capacity at constant volume. This can sometimes cause confusion in the obtained results. See also \cite{thermo}.\\

The following relations will be proven by considering a divariant substance:
\begin{equation}\label{Relation_5}
    d\,e=C_v\,d\,T+\frac{\gamma\,\Gamma\,p-(\gamma-1)\,\rho\,c^2}{\gamma\,\Gamma\,\rho^2}d\,\rho\,,
\end{equation}
\begin{equation}\label{Relation_6}
   d\,h=\gamma\,C_v\,d\,T+\frac{\Gamma-\gamma+1}{\Gamma\,\rho}\,d\,p\,,\quad h=e+\frac{p}{\rho}\,,
\end{equation}
\begin{equation}\label{Relation_7}
    d\,p=c^2\,d\,\rho+\rho\,\Gamma\,T\,\,d\,s=\frac{c^2}{\gamma}\,d\,\rho+\rho\,\Gamma\,C_v\,d\,T\,.
\end{equation}
Indeed, starting from:
\begin{equation}\label{Relation_1}
    d\,e=C_v\,d\,T+\left(\beta+\frac{p}{\rho^2}\right)d\,\rho\,,
\end{equation}
\begin{equation}\label{Relation_2}
   d\,h=\gamma\,C_v\,d\,T+\left(\alpha+\frac{1}{\rho}\right)d\,p\,,\quad h=e+\frac{p}{\rho}\,,
\end{equation}
\begin{equation}\label{Relation_3}
    d\,p=c^2\,d\,\rho+\rho\,\Gamma\,T\,\,d\,s=c_T^2\,d\,\rho+\epsilon\,d\,T\,,
\end{equation}
where $\alpha$\,, $\beta$\,, \ldots\, can be seen as partial derivatives {\it e.g.} $\epsilon\equiv\frac{\partial p}{\partial T}\big|_\rho$ and where the pressure $p$ for this divariant substance is seen as a function of the two independent thermodynamic variables $\rho$ and $T$. Some of the coefficients have already been identified in order to be consistent with Table \ref{tab4}, \eqref{CT} and \eqref{Cv_etc}.
It is elementary to prove the following result using these two variables.
\begin{prop}\label{prop30}
It follows from Gibbs relation:
\begin{equation}\label{Gibbs}
    T\,d\,s=d\,e-\frac{p}{\rho^2}d\,\rho\,,
\end{equation}
that
\begin{multline}\label{Relation_4}
    \epsilon=\Gamma\,\rho\,C_v \,,\quad \alpha=-\frac{\gamma-1}{\Gamma\,\rho}\,,\quad c_T^2=\frac{c^2}{\gamma}\,,\quad \beta=-\frac{(\gamma-1)\,c^2}{\gamma\,\Gamma\,\rho}\,,\\\Gamma^2\,C_v\,T=\frac{\gamma-1}{\gamma}c^2\,.
\end{multline}
\end{prop}
Hence, by combining identities \eqref{Relation_5} to \eqref{Relation_7} it is found that:
\begin{equation}\label{Relation_8}
d\,p=c^2\,d\,\rho+\rho\,\Gamma\,T\,d\,s\,,\quad d\,T=\frac{\Gamma\,T}{\rho}\,d\,\rho+\frac{T}{C_v}\,d\,s\,,
\end{equation}
which is needed for the derivation of \eqref{systeme_NST_L}.

\newpage


\bibliographystyle{model1-num-names}



\newpage
\tableofcontents
\end{document}